\documentclass[a4paper]{article}

\usepackage[english]{babel}
\usepackage[utf8x]{inputenc}
\usepackage[T1]{fontenc}
\usepackage[a4paper,top=3cm,bottom=2cm,left=3cm,right=3cm,marginparwidth=1.75cm]{geometry}
\usepackage{amsmath}
\usepackage{graphicx}
\usepackage[colorinlistoftodos]{todonotes}
\usepackage{listings}
\usepackage{color} %red, green, blue, yellow, cyan, magenta, black, white
\definecolor{mygreen}{RGB}{28,172,0} % color values Red, Green, Blue
\definecolor{mylilas}{RGB}{170,55,241}
\usepackage{changepage}
\usepackage{graphicx}
\usepackage[small]{caption}
\usepackage{subcaption}
\usepackage[makeroom]{cancel}
\usepackage{listings}
\usepackage{amsthm}
\usepackage{amsfonts}
\usepackage{float}
\usepackage{listings}
\usepackage{framed}
\usepackage{cite}
\usepackage{algorithm}
\usepackage{pdfpages}
\usepackage{mathtools}
\usepackage{bm}
\usepackage{authblk}
\usepackage{multirow}
\usepackage{authblk}
\usepackage{tabularx}
\usepackage{pbox}

\usepackage{todonotes}

\DeclareMathOperator{\arctantwo}{arctan2}
\newcommand{\norm}[1]{\left\lVert#1\right\rVert}

\begin{document}
	
\title{Efficiency and Parameter Selection of a micro-macro Markov chain Monte Carlo method for molecular dynamics}

\author[1]{Hannes Vandecasteele}
%\ead{hannes.vandecasteele@cs.kuleuven.be}
\author[1]{Giovanni Samaey}
%\ead{giovanni.samaey@cs.kuleuven.be}
\affil[1]{KU Leuven, Department of Computer Science, NUMA Section, Celestijnenlaan 200A box 2402, 3001 Leuven, Belgium}
\date{\today}
	
\maketitle

\begin{abstract}
 We recently introduced a mM-MCMC scheme that is able to accelerate the sampling of Gibbs distributions when there is a time-scale separation between the complete molecular dynamics and the slow dynamics of a low dimensional reaction coordinate. The mM-MCMC Markov chain works in three steps: 1) compute the reaction coordinate value associated to the current molecular state; 2) generate a new macroscopic proposal using some approximate macroscopic distribution; 3) reconstruct a molecular configuration that is consistent with the newly sampled macroscopic value. There are a number of method parameters that impact the efficiency of the mM-MCMC method. On the macroscopic level, the proposal- and approximate macroscopic distributions are important, while on the microscopic level the reconstruction distribution is of significant importance.  In this manuscript, we will investigate the impact of these parameters on the efficiency of the mM-MCMC method on three molecules: a simple three-atom molecule, butane and alanine-dipeptide.
\end{abstract}

\section{Introduction}
The main objective of this manuscript is to study the numerical effiency of a newly proposed micro-macro Markov chain Monte Carlo (mM-MCMC) method~\cite{vandecasteele2022}. We designed this method with the goal of accelerating the sampling of molecular dynamics in mind. Many molecular systems typically exist of a large number of interacting atoms, whose motion is described by a complex potential energy
\begin{equation} \label{eq:potenergy}
V\left(x\right) = V_{\text{bonds}}(x)+V_{\text{angles}}(x) + V_{\text{torsions}}(x) + V_{\text{electrostatic}}(x) + V_{\text{LJ}}(x).
\end{equation}
In this formula, $x\in \mathbb{R}^{3d}$ is the vector containing the tree-dimensional positions of the atoms. When there are no chemical reactions, the potential energy can usually be composed as a sum of local and non-local potentials. For every covalent bond and angle between three connected atoms there is a quadratic potential, $V_{\text{bonds}}$ and $V_{\text{angles}}$ respectively, that describes vibrations of the atoms around a equilibrium value. Similarly, certain torsion angles between four connected atoms are described by a periodic potential $V_{\text{torsions}}$. Finally, there are two important non-local potentials, one for electrostatic interactions, $V_{\text{electrostatic}}$, and another for Vanderwaals interactions. The latter is mostly modelled by a Lennard-Jones potential $V_{\text{LJ}}$. With this potential energy, we can define the time-invariant Gibbs distribution
\begin{equation} \label{eq:gibbs}
\mu(x) = Z_V^{-1} \exp\left(-\beta V(x)\right) \ \ \  \forall \ x \in \mathbb{R}^{3d},
\end{equation} 
that describes the equilibrium behaviour of the molecule. The constant $\beta$ is the inverse temperature, and $Z_V$ is the normalization constant.	

There are two serious limitations when sampling from this invariant distribution. First, there can be a time-scale separation $\varepsilon$ between the macroscopic, global conformation of the molecule and the microscopic vibrations of the individual atoms. In this situation, the potential energy decomposes as a sum of slow and fast terms
\begin{equation*} \label{eq:Vmultiscale}
V(x) = V_{s}(x) + \frac{1}{\varepsilon} V_{f}(x),
\end{equation*}
where $V_s$ is the potential energy of the slow, macroscopic conformation, and $V_f$ is the potential energy of the remaining fast degrees of freedom. Both are independent from $\varepsilon$. The existence of the time-scale separation calls for the use of small time steps on the order of $\varepsilon$ to keep the simulation stable. However, these small time steps prohibit sampling all possible global structures the molecular can have. Second, there is also the problem of dimensionality. It has been shown that the acceptance rate of some Markov chain Monte Carlo decreases as $d^{-1/3}$, where $d$ is the dimension of the system. We will ignore the issue of large dimensionality in this paper.

Instead, we will focus on cases where the macroscopic system is defined by the use of reaction coordinates. A reaction coordinate is a differentiable function $\xi$ that maps a microscopic molecule to some macroscopic description, i.e., 
\begin{equation} \label{eq:rc}
\xi : \mathbb{R}^d \to \mathbb{R}^n : x \mapsto z,
\end{equation}
with $n \ll d$. Given a certain reaction coordinate, we can derive a potential energy and invariant Gibbs distribution for it. Indeed, the potential energy, also called the free energy~\cite{le2012mathematical}, of $\xi $ is given by
\begin{equation} \label{eq:freeenergy}
A(z) = -\frac{1}{\beta} \int_{\mathbb{R}^d} Z_V^{-1} \exp\left(-\beta V(x)\right) \delta_{\xi(x)-z}(dx).
\end{equation}
Note that the Dirac delta function only considers $x$ values that have reaction coordinate value $z$. Then, similarly to~\eqref{eq:gibbs}, the reaction coordinate values follow a time-invariant Gibbs distribution of the form
\begin{equation} \label{eq:mu0}
\mu_0(z)  = Z_A^{-1} \exp\left(-\beta A(z)\right).
\end{equation}
With this free energy, we can rewrite the total microscopic potential energy~\eqref{eq:Vmultiscale} as
\begin{equation}
V(x) = A\left(\xi(x)\right) + \frac{1}{\varepsilon} V_f(x).
\end{equation}
Given both the macroscopic and microscopic Gibbs distributions, we can define the time-invariant reconstruction distribution from macro to micro~\cite{le2012mathematical}. This reconstruction distribution is defined on the level set $\Sigma(z) = \{ x \in \mathbb{R}^d\ | \ \xi(x) = z \}$. This distribution is defined as
\begin{equation} \label{eq:exactnu}
\nu(x | z) = \frac{\mu(x)}{\mu_0(z)}, \ x \in \Sigma(z).
\end{equation}

The time-scale separation problem is far from new. The past decades have seen many creative methods to overcome this problem and reduce its cost. We give an overview. The Adaptive Biasing methods~\cite{dickson2010free,comer2015adaptive} flatten the energy in the direction of the reaction coordinate. There are two variants of this method; adaptive biasing potential (ABP)~\cite{berg1991multicanonical,wang2001efficient, laio2002escaping,marsili2006self,zheng2008random} and adaptive biasing force (ABF)~\cite{darve2001calculating,darve2008adaptive}. Both achieve improved sampling by updating a running approximation to the free energy (ABP) or mean force (ABF) and subtracting it from the microscopic potential/force. We also like to mention the recently developped Multilevel Markov Chain Monte Carlo method. There are two instantiations that reduce the variance of MCMC by coupling levels. Dodwell et. al.~\cite{dodwell2019multilevel}, couple two levels by using the delayed-acceptance method~\cite{christen2005markov}. That is, particle proposals are generated on the coarser level, and then only accepted or rejected on the finer level. A second approach to Multilevel MCMC is to use the same random numbers for proposals on both levels~\cite{jasra2018markov}. It is unclear which scheme can lead to the largest time/variance reduction. Finally, The Coupled MCMC method for lattice systems~\cite{kalligiannaki2012coupled} also accelerates sampling on a fine lattice by coupling it with a coarser lattice, like above MLMCMC methods. However, the coarser lattice does not live in the same state space as the fine lattice. In this method, going back to a finer level is done by sampling from a reconstruction distribution that is constrained to the coarser state.

In previous work, we created a micro-macro MCMC method~\cite{vandecasteele2022} (mM-MCMC) similar to the coupled MCMC method~\cite{kalligiannaki2012coupled}. With the current microscopic sample, it uses three steps to generate a new microscopic sample: i) compute its reaction coordinate value; ii) Sample a new macroscopic reaction coordinate value; iii) Build a microscopic sample by sampling from a (constrained) reconstruction distribution. The mM-MCMC method has the same structure as the delayed-acceptance, with the difference that the macroscopic level has a lower dimension. We also proved convergence and demonstrated that our method can realize substantial efficiency gains. There are, nonetheless, many parameters that influence these gains. We study the optimal choice of these parameters in this manuscript. 

\paragraph{Outline of the manuscript}
We start by introducing the mM-MCMC method in Section~\ref{sec:2}, followed by the criterion used to measure the performance of the mM-MCMC method in Section~\ref{sec:effgain}. In Section~\ref{sec:threeatom}, we study the three-atom molecule, where there is a single parameter $\varepsilon$ that determines the time-scale separation. We investigate the performace of both mM-MCMC schemes on this example and look at the effect of the method parameters on the efficiency gain. Next, in Section~\ref{subsec:butane}, we study the parameters of direct and indirect reconstruction on the more complex molecule butane, and end with a comparison between both schemes.

\section{The micro-macro Markov chain Monte Carlo method} \label{sec:2}
We now describe the mM-MCMC method in more detail. In its most general form, the mM-MCMC method consists of three steps, i) Large macroscopic proposal move followed by a macroscopic acceptance/rejection step, ii) Reconstruction of a new microscopic sample that is consistent with the new reaction coordinate value, iii) microscopic acceptance/rejection step to make sure the new microscopic sample is distrbuted according to~\eqref{eq:gibbs}.

In previous work~\cite{vandecasteele2022}, we introduced two variants of the mM-MCMC method: with direct reconstruction and with indirect reconstruction. As the name suggests, the two methods differ only in the reconstruction step. The mM-MCMC method with direct reconstruction always generates a new microscopic sample on the sub-manifold of constant reaction coordinate. The mM-MCMC method with indirect reconstruction, on the other hand, generates a microscopic sample by pulling the previous microscopic sample towards, but not necessarily on, the sub-manifold of constant reaction coordinate. The direct reconstruction variant is especially useful when the reaction coordinate is simple and $\Sigma(z)$ has a regular shape. In other cases, the indirect reconstruction is easier to use.

\subsection{mM-MCMC with direct reconstruction}
Let us now fill in the three steps from above in case of mM-MCMC with direct reconstruction.

\subsubsection{Macroscopic proposal move} \label{subsubsec:macropropdirect}
Given the current microscopic sample $x$ and its reaction coordinate $z$, we generate a new reaction coordinate value $z'$ using a macroscopic proposal distribution $q(z, z')$. This conditional distribution is usually based on a macroscopic time stepping scheme. Examples include a gradient descent step (MALA), a simple discrete Brownian Motion, or more complex time stepping schemes. As a consequence, $q(z, z')$ is usually a Gaussian distribution in $z'$ with a mean related to $z$.

\subsubsection{Macroscopic acceptance}
Before proceeding to the reconstruction step, we need to make sure that $z'$ is distributed according to some approximate macroscopic distribution $\bar{\mu}_0$. This approximate distribution can be any distribution chosen by the user since the exact time-invariant distribution~\eqref{eq:mu0} is not readily availble for complex molecules. Therefore, we accept $z'$ with probability
\begin{equation} \label{eq:alphacg}
\alpha_{CG}(z, z') = \min\left\{1, \frac{\bar{\mu}_0(z') \ q_0(z', z)}{\bar{\mu}_0(z) \ q_0(z, z')} \right\}.
\end{equation}
This is the standard Metropolis-Hastings acceptance probability~\cite{metropolis1953equation,hastings1970monte}.

\subsubsection{Reconstruction}
Given the new macroscopic value $z'$, we reconstruct a microscopic sample $x'$ by drawing a sample from the reconstruction distribution $\bar{\nu}(\cdot | z')$ on the sub-manifold $\Sigma(z')$. This distribution can be anyting specified beforehand by the user. 

\subsubsection{Microscopic acceptance}
Finally, to ensure that the microscopic sample is indeed distributed according to~\eqref{eq:gibbs}, we accept $x'$ with probability
\begin{equation*} \label{eq:alpha_f_direct}
\alpha(x, x') = \min\left\{1, \frac{\mu(x') \bar{\mu}_0(z) \bar{\nu}(x|z)}{\mu(x) \bar{\mu}_0(z') \bar{\nu}(x'|z')}  \right\}.
\end{equation*}
With this microscopic acceptance rate, it can be shown that the mM-MCMC method indeed has~\eqref{eq:gibbs} as its invariant distribution and that it is ergodic.

\subsection{mM-MCMC with indirect reconstruction}
For technical reasons, the mM-MCMC method with indirect reconstruction is defined on the extended state space $(x, z) \in \mathbb{R}^d \times \mathbb{R}^n$ of microscopic and reaction coordinate values. See~\cite{vandecasteele2022} for more details. Therefore, there is no need to compute the reaction coordinate value of the current microscopic sample $x$. Let us now explain the different steps of the algorithm.

\subsubsection{Macroscopic proposal move \& acceptance}
In this step, we again generate a new reaction coordinate value $z'$ with the transition distribution $q_0(z, z')$. Afterwards, we accept $z'$ with probability $\alpha_{CG}$~\eqref{eq:alphacg}. This step is identical to section~\ref{subsubsec:macropropdirect}.

\subsubsection{Reconstruction}
Instead of generating a microscopic sample $x'$ on the manifold $\Sigma(z')$, we sample $x'$ from the \textit{indirect reconstruction distribution}
\begin{equation}
\nu_{\lambda}(x', z') = N_{\lambda}^{-1}(z') \left(\frac{\lambda \beta}{2\pi}\right)^{n/2} \exp\left(-\beta V(x')\right) \exp\left(-\frac{\lambda \beta \norm{\xi(x') - z'}^2}{2}\right).
\end{equation}
The normalization constant $N_{\lambda}(z')$ is defined as
\begin{equation}
N_{\lambda}(z') = Z_V\int_{\mathbb{R}^n} \left(\frac{\lambda \beta}{2\pi}\right)^{n/2} \exp\left(-\frac{\lambda \beta}{2}(u-z')^2\right) d\mu_0(u).
\end{equation}
Note that this indirect reconstruction distribution is centered around $\Sigma(z)$ and is almost zero elsewhere.

In this manuscript, we will take a sample from this indirect reconstruction distribution by running a short Markov chain Monte Carlo process. Specifically, we will take $K$ time steps of the MALA algorithm with time step $\delta t$:
\begin{equation}
x_{k+1} = x_k - \delta t \nabla V(x_k) - \delta t \lambda (\xi(x_k) - z') \nabla \xi(x_k) + \sqrt{2\delta t \beta^{-1}} \eta_k,  \ \eta_k \sim \mathcal{N}(0, 1),  \  k = 0,\dots, K-1,
\end{equation}
and set $x' = x_K$. This way, if $K$ is large enough, the microscopic particle $x'$is distributed according to $\nu_{\lambda}(\cdot, z')$. The parameter $\lambda$ determines the speed at which $\xi(x_k)$ converges to $z'$. The higher $\lambda$, the faster the convergence, making it the most important parameter in the mM-MCMC method with indirect reconstruction.

\subsubsection{Microscopic acceptance}
Finally, we must ensure that $x'$ is indeed a sample from the Gibbs distribution~\eqref{eq:gibbs}. We ensure this by accepting $x'$ with probability
\begin{equation*}
\alpha(x, x') = \min\left\{1, \frac{\bar{\mu}_0(z) N_{\lambda}(z')}{\bar{\mu}_0(z') N_{\lambda}(z)} \right\}.
\end{equation*}
Note that this acceptance rate only depends on the reaction coordinate value $z$ and $z'$, in contrast to~\eqref{eq:alpha_f_direct}. This form is because of the explicit form of the indirect reconstruction distribution.

\section{Efficiency Gain Criterion} \label{sec:effgain}
To assess the total efficiency gain of the mM-MCMC method, we will always compare its sampling result with the microscopic Metropolis Adjusted Langevin (MALA) method~\cite{xifara2014langevin}. We briefly explain this method for self-containement in Section~\ref{subsec:mala}. Afterwards, we introduce the efficiency criterion in Section~\ref{subsec:effgain}.

\subsection{The Metropolis adjusted Langevin method} \label{subsec:mala}
MALA is a well known method to sample the time-invariang Gibbs distribution~\eqref{eq:gibbs}. Given the current microscopic sample $x_n \in \mathbb{R}^d$, at time $t_n = n \delta t$, it generates the next sample as
\begin{equation} \label{eq:mala}
x_{n+1} = x_n - \delta t \nabla V(x_n) + \sqrt{2dt\beta} \eta_n,
\end{equation}
with $\eta_n$ distrbuted according to the $d-$dimensional normal distribution with mean zero and unit covariance. We then accept $x_{n+1}$ with Metropolis-Hastings probability
\begin{equation*}
\alpha(x_n, x_{n+1}) = \min\left\{1, \frac{\mu(x_{n+1}) q(x_{n+1}, x_n)}{\mu(x_n) q(x_n, x_{n+1})}  \right\}.
\end{equation*}
In the above formula, $q(x_n, x_{n+1})$ is the transition probability distribution associated to~\eqref{eq:mala}.

\subsection{Efficiency gain} \label{subsec:effgain}
Suppose we are interested in estimating the expected value of some functional $F: \mathbb{R}^d \to \mathbb{R}$ with respect to the Gibbs measure, 
\begin{equation} \label{eq:F}
\mathbb{E}[F] = \int_{\mathbb{R}^d} F(x) d\mu(x).
\end{equation}
By drawing random samples from $\mu$, through a Markov chain Monte Carlo method, we can obtain an estimate for $\mathbb{E}[F]$. To assess the accuracy of the MCMC method, one can perform $M$ independent runs, each with estimated value $F_i$ for~\eqref{eq:F}, and compute the Mean Squared Error (MSE).
\begin{equation}
\text{MSE} = \frac{1}{M} \sum_{i=1}^M \left( F_i -\mathbb{E}[F] \right)^2.
\end{equation}

In all coming experiments, we will compare the performace of the mM-MCMC method to the microscopic MALA method. Therefore, we define the efficiency gain of mM-MCMC over MALA as
\begin{equation} \label{eq:effgain}
\frac{\text{MSE}_{\text{MALA}}}{\text{MSE}_{\text{mM-MCMC}}} \frac{T_{\text{MALA}}}{T_{\text{mM-MCMC}}},
\end{equation}
where $T_{\text{MALA}}$ and $T_{\text{mM-MCMC}}$ are the respective CPU times of the MALA and mM-MCMC methods.

\section{The three-atom molecule} \label{sec:threeatom}
\begin{figure}
	\centering
	\includegraphics[width=0.3\linewidth]{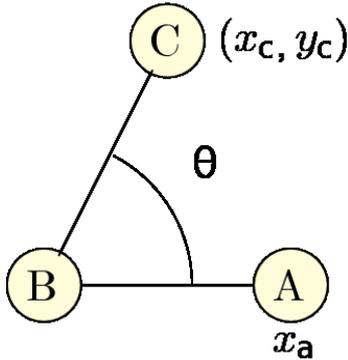}
	\caption{The three-atom molecule. Atom $A$ is constraint to the $x$-axis with $x$-coordinate $x_a$, atom $B$ is fixed at the origin of the plane and atom $C$ lies om the two-dimensional plane with Cartesian coordinates $(x_c, y_c)$.}
	\label{fig:triatommolecule2}
\end{figure}

In this section, we consider the mM-MCMC algorithm on a simple, academic, three-atom molecule, as first introduced in~\cite{legoll2012some}. The three-atom molecule has a central atom $B$, that we fix at the origin of the two-dimensional plane, and two outer atoms, $A$ and $C$. To fix the superfluous degrees of freedom, we constrain atom $A$ to the $x-$axis, while $C$ can move freely in the plane. The three-atom molecule is depicted on Figure~\ref{fig:triatommolecule2}. 

The potential energy for the three-atom system consists of three terms,
\begin{equation} \label{eq:trheeatompotential}
V(x_a, x_c, y_c) = \frac{1}{2\varepsilon} \ (x_a-1)^2 + \frac{1}{2\varepsilon} \ (r_c-1)^2 + \frac{208}{2}\left(\left(\theta-\frac{\pi}{2}\right)^2 - 0.3838^2\right)^2,
\end{equation}
where $x_a$ is the $x-$coordinate of atom $A$ and $(x_c, y_c)$ are the Cartesian coordinates of atom $C$. The bond length $r_c$ between atoms $B$ and $C$ and the angle $\theta$ between atoms $A$, $B$ and $C$ are defined as
\begin{equation*}
\begin{aligned}
r_c &= \sqrt{x_c^2 + y_c^2} \\ 
\theta &= \arctantwo (y_c, x_c).
\end{aligned}
\end{equation*}
The first term in~\eqref{eq:trheeatompotential} describes the vibrational potential energy of the bond between atoms $A$ and $B$, with equilibrium length $1$. Similarly, the second term describes the vibrational energy of the bond between atoms $B$ and $C$ with bond length $r_c$. Finally, the third term determines the potential energy of the angle $\theta$ between the two outer atoms, which has an interesting bimodal behaviour. The distribution of $\theta$ has two peaks, one at $\frac{\pi}{2} -0.3838$ and another at $\frac{\pi}{2} + 0.3838$. 

The reaction coordinate that we consider in this section is the angle $\theta$, i.e.,
\begin{equation}
\xi(x_a, x_c, y_c) = \theta(x_a, x_c, y_c),
\end{equation}
since this variable is the slow component of the three-atom molecule. Additionally, the angle $\theta$ is also independent of the time-scale separation, given by $\varepsilon$.

Finally, from the potential energy~\eqref{eq:trheeatompotential}, it can be easily seen that the free energy of the reaction coordinate $\theta$ reads
\begin{equation} \label{eq:threeatomA}
A(z) = \frac{208}{2}\left(\left(z-\frac{\pi}{2}\right)^2 - 0.3838^2\right)^2.
\end{equation}

\paragraph{\textbf{Outline of this section}}
In previous work~\cite{vandecasteele2022}, we have demonstrated that both mM-MCMC methods with direct and indirect reconstruction can sample the invariant distribution~\eqref{eq:trheeatompotential} correctly. We continue this line of experiments in this manuscript. Specifically, we include the following numerical results.

We start with all experiments relating to the direct reconstruction scheme in Section~\ref{subsec:threeatomdirect}. First, we show in Section~\ref{subsubsec:comparsionmacro} that the mM-MCMC method with direct reconstruction is able to sample the free energy correctly, for three separate choices for the approximate macroscopic distribution $\bar{\mu}_0$. Then, in Section~\ref{subsubsec:Aq} we study the effect of $\bar{\mu}_0$ and $q_0$ on the efficiency gain of the mM-MCMC method. This is followed by an investigation of $\bar{\mu}_0$ and $\bar{\nu}$ on its efficiency gain in Section~\ref{subsubsec:Anu}. Finally, we study the efficiency gain of mM-MCMC with direct reconstruction as a function of $\varepsilon$, Section~\ref{subsubsec:effgainepsdirect}.

Afterwards, we study the mM-MCMC method with indirect reconstruction in Section~\ref{subsec:threeatomindirect}. First, we see what values of $\lambda$ are optimal for reconstruction in case of the three-atom molecule. Then, with $\lambda$ we study the efficiency gain as a function of $\varepsilon$, respectively in Sections~\ref{subsubsec:threeatomlambda} and~\ref{subsubsec:triatomindriecteps}. Finally, we compare the mM-MCMC methods with direct and indirect reconstruction with each other in Section~\ref{subsubsec:comparisondirectindirect}.

\subsection{Experiments with direct reconstruction} \label{subsec:threeatomdirect}

\subsubsection{Comparison with the macroscopic reaction coordinate sampler} \label{subsubsec:comparsionmacro}
In the second numerical illustration, we compare the mM-MCMC method to MCMC samplers at the macroscopic level. That is, if the approximate macroscopic distribution $\bar{\mu}_0$ differs significantly from the exact marginal distribution of the reaction coordinate $\mu_0$~\eqref{eq:threeatomA}, any macroscopic MCMC sampler would render a poor fit on $\mu_0$. However, the mM-MCMC method does not readily accept reaction coordinate values sampled at the macroscopic level, but performs a second accept/reject step at the microscopic level. This second accept/reject step guarantees that the reconstructed microscopic samples are indeed samples from the microscopic Gibbs measure $\mu$~\eqref{eq:gibbs}. As a consequence, the reaction coordinate values sampled with mM-MCMC will sample the marginal distribution $\mu_0$ exactly.

\begin{figure}
	\centering
	\includegraphics[width=0.55\linewidth]{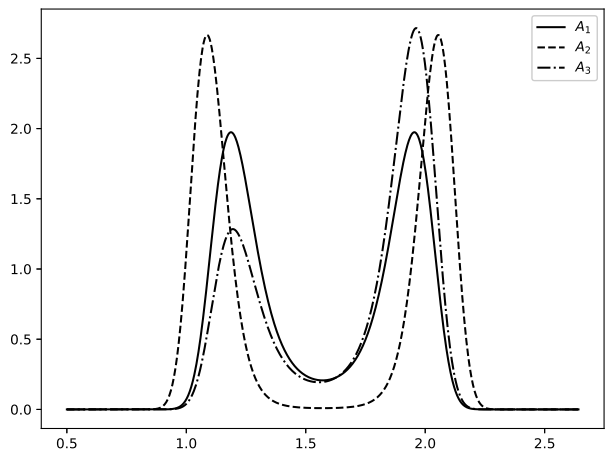}
	\caption{The three approximate free energy distributions as defined in equation~\eqref{eq:freeenergies_2}. The blue curve is the exact free energy of the three-atom molecule~\eqref{eq:trheeatompotential}, the red curve is obtained by shifting the local maxima of the exact distribution by $0.1$ by the left and the right, while the final distribution is the exact expression perturbed by the cosine function.}
	\label{fig:freeenergy}
\end{figure}

\paragraph{Experimental setup} We verify this claim by defining two approximate free energy functions $\bar{A}$ that define two approximate macroscopic distributions $\bar{\mu}_0$. These approximate free energy functions are
\begin{align} \label{eq:freeenergies_2}
\bar{A}^1(\theta) &= \frac{208}{2}\left(\left(\theta-\frac{\pi}{2}\right)^2 - 0.4838^2\right)^2, \\
\bar{A}^2(\theta) &= \frac{208}{2}\left(\left(\theta-\frac{\pi}{2}\right)^2 - 0.3838^2\right)^2 + \cos(\theta). \nonumber
\end{align}
The first expression for the approximate free energy is obtained by perturbing the two peaks of the exact free energy~\eqref{eq:freeenergy} by a distance $0.1$ of radians. The second formula, $\bar{A}^3$ perturbs the exact free energy expression by adding the cosine function to $A$, resulting in a large perturbation on the amplitude of the associated macroscopic invariant distribution. The effect of the cosine perturbation is that the hight of the left peak in the macroscopic distribution of $\theta$ is decreased, while the height of the right peak is increased. The three Gibbs distributions associated to the correct free energy and the two approximate free energies defined above are shown in Figure~\ref{fig:freeenergy}.

In the following experiment, we investigate the numerical sampling results of the mM-MCMC method where we use each of the above defined approximate macroscopic distributions. On the macroscopic level, we use the MALA method~\eqref{eq:mala} with time step $\Delta t=0.01$ to generate reaction coordinate values according to each of the approximate macroscopic distributions. For reconstruction, we use the time-invariant reconstruction distribution $\nu$~\eqref{eq:exactnu}. We also choose the time-scale parameter $\varepsilon=10^{-6}$ and sample $N=10^6$ microscopic samples. 

\begin{figure}
	\centering
	\begin{subfigure}[b]{0.5\textwidth}
		\centering
		\includegraphics[width=1.\linewidth]{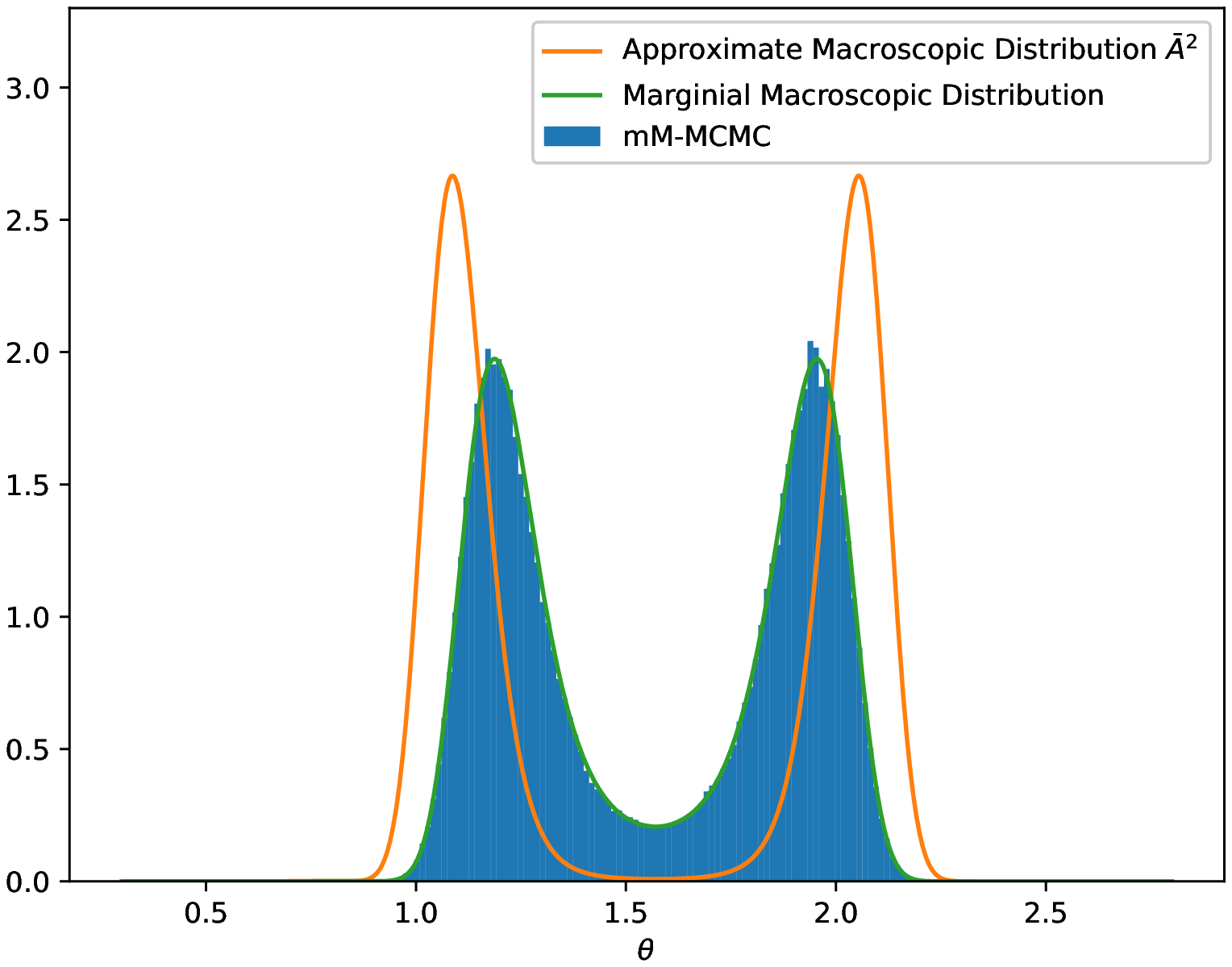}
	\end{subfigure}%
	\begin{subfigure}[b]{0.5\textwidth}
		\centering
		\includegraphics[width=1.\linewidth]{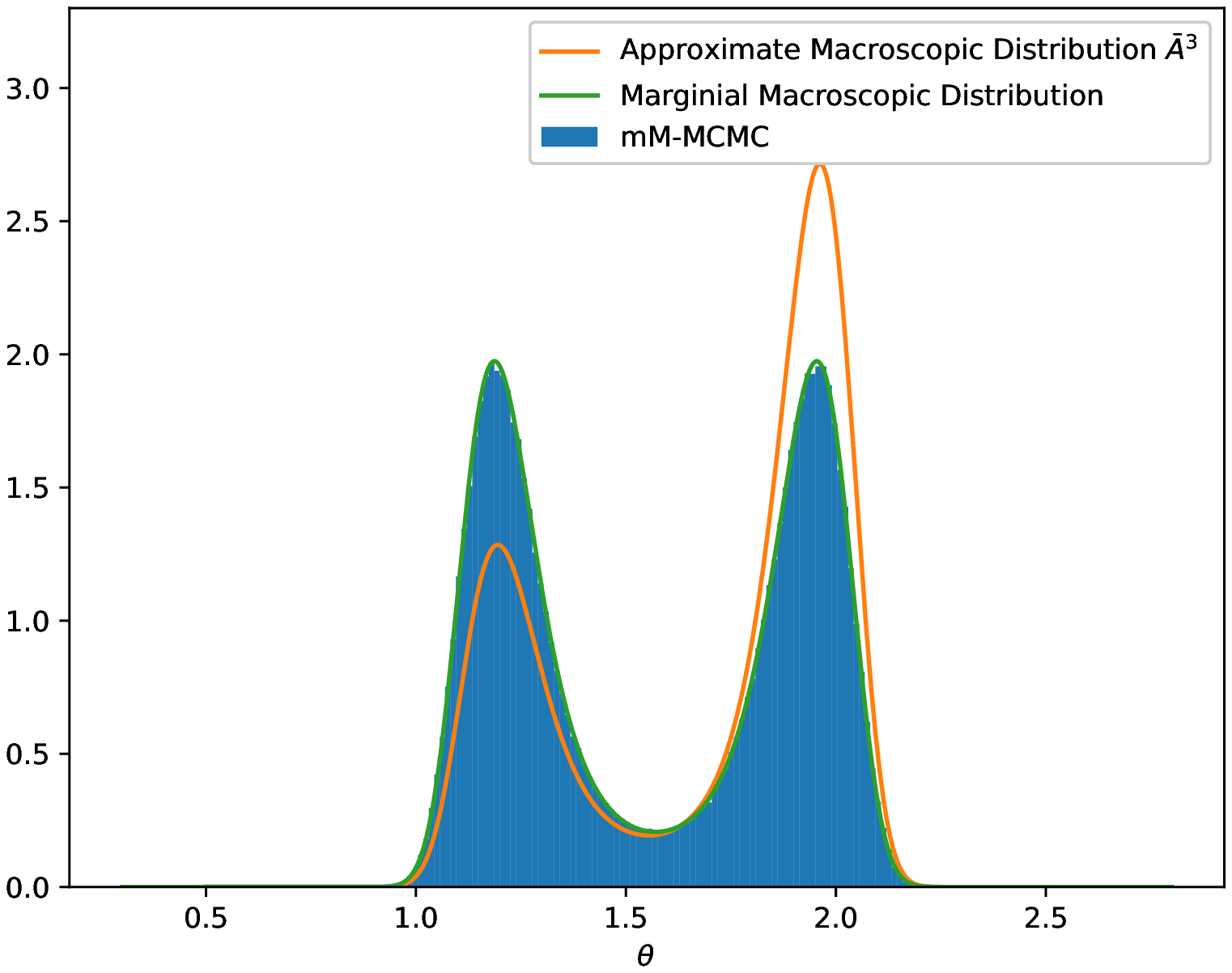}
	\end{subfigure}
	\caption{Histogram of the mM-MCMC method (blue) for two different approximate macroscopic distributions: $\bar{A}^2$ (left) and $\bar{A}^3$ (right). The orange line depicts the approximate macroscopic distribution used in each experiment, and the green line represents the marginal distribution of the reaction coordinate values. We can see that the mM-MCMC is able to overcome a bad approximate macroscopic distribution and samples the marginal distribution correctly.}
	\label{fig:macro_approx}
\end{figure}

\paragraph{Numerical results} The mM-MCMC method is in both cases able to recover the marginal distribution of the reaction coordinates. A simple macroscopic sampler would return samples that differ significantly from the correct macroscopic marginal distribution. This result thus shows the need for a microscopic correction to an inexact macroscopic sampler, and the merit of the mM-MCMC method.

\subsubsection{Impact of $\bar{A}$ and $q_0$ on the efficiency gain} \label{subsubsec:Aq}
\paragraph{\textbf{Experimental setup}}
Having visually shown that the mM-MCMC can correct for an inaccurate approximate macroscopic distribution, we numerically investigate the impact of the choice of approximate macroscopic distribution $\bar{\mu}_0$ and the choice of macroscopic transition distribution $q_0$ on the efficiency gain of mM-MCMC.

For the approximate macroscopic distribution, we use the exact free energy function~\eqref{eq:threeatomA} and the two approximate macroscopic distribution that we defined in~\eqref{eq:freeenergies_2}. The two choices for the macroscopic proposal distribution $q_0$ are based on the following two stochastic dynamical systems for the reaction coordinate,
%\begin{equation*}
\begin{align} \label{eq:coarseproposalmoves}
q_0^1  &:d\theta = -\nabla A_i(\theta)dt + \sqrt{2\beta^{-1}} dW \nonumber \\ %%%%%%%%
q_0^2 &: d\theta = \sqrt{2\beta^{-1}} dW.
\end{align}
%\end{equation*}
The first stochastic differential equation is the overdamped Langevin dynamics for each of the approximate free energy functions (MALA), while the latter equation is a simple Brownian motion in the reaction coordinate space.

In the following experiment, we run the mM-MCMC algorithm with each of these six combinations for the approximate macroscopic invariant distribution $\bar{\mu}_0 \propto \exp\left(-\beta\bar{A}\right)$ and the macroscopic transition distribution $q_0$ for $N=10^6$ sampling steps and with a macroscopic time step $\Delta t = 0.01$. We also take the reconstruction distribution $\nu$ in each of these experiments. After each simulation with these combinations, we compute the total macroscopic acceptance rate, the microscopic acceptance rate after reconstruction, the average runtime, the numerical variance on the estimated mean of $\theta$ and the total efficiency gain of of mM-MCMC over the MALA algorithm, as explained in Section~\ref{subsec:effgain}. The microscopic time step for the MALA method is $\delta t = \varepsilon$ and for we choose $\beta=1$ for the inverse temperature. For a statistically good comparison, we average the results over $100$ independent runs. The numerical efficiency gains are depicted in Tables~\ref{tab:threeatomgain1e-04} and~\ref{tab:threeatomgain1e-06} for $\varepsilon=10^{-4}$ and $\varepsilon=10^{-6}$ respectively.

\begin{table}[h]
	\centering
	\begin{tabular}{c|c|c|c|c|c}
		\centering
		Parameters & \pbox{15cm}{Macroscopic \\ acceptance rate} & \pbox{15cm}{Microscopic \\ acceptance rate} & \pbox{15cm}{Runtime \\ gain} & \pbox{15cm}{Variance \\ gain} & \pbox{15cm}{Total \\ efficiency gain} \\
		\hline
		Langevin, $\bar{A}^1$ & 0.749932            &           1            &      2.45692     &      85.3266       &             209.64 \\
		Langevin, $\bar{A}^2$ &                      0.730384              &         0.432508      &     2.62306  &        28.4122          &           74.527 \\
		Langevin, $\bar{A}^3$    &                   0.749653         &              0.950238     &      1.95663       &    99.7186            &       195.112 \\
		Brownian, $\bar{A}^1$    &                   0.645188           &            1                 & 3.03108        &      85.1586               &          258.122 \\
		Brownian, $\bar{A}^2$     &                 0.61375             &           0.597058       &    3.28045 &         35.6794         &             117.044 \\
		Brownian, $\bar{A}^3$        &               0.645654            &           0.959794       &    2.72728     &       81.3405        &           221.838 \\
	\end{tabular}
	\caption{A summary of different statistics of the mM-MCMC method with $\varepsilon=10^{-4}$ for six combinations of the (approximate) macroscopic invariant distribution and macroscopic proposal moves. We record the macroscopic and microscopic acceptance rates of mM-MCMC, the runtime and variance gains of mM-MCMC over the microscopic MALA method and the total efficiency gain.}
	\label{tab:threeatomgain1e-04}
\end{table}

\begin{table}[h]
	\centering
	\begin{tabular}{c|c|c|c|c|c}
		\centering
		Parameters & \pbox{15cm}{Macroscopic \\ acceptance rate} & \pbox{15cm}{Microscopic \\ acceptance rate} & \pbox{15cm}{Runtime \\ gain} & \pbox{15cm}{Variance \\ gain} & \pbox{15cm}{Total \\ efficiency gain} \\
		\hline
		Langevin, $\bar{A}^1$         &              0.749906        &               1           &        2.50343       &     3297.65           &       8255.44 \\
		Langevin, $\bar{A}^2$           &             0.730538              &         0.43242    &        2.64621    &    933.64    &      2470.61 \\
		Langevin, $\bar{A}^3$           &           0.749596               &        0.950308    &       1.99301     &  2461.08   &           4904.96 \\
		Brownian, $\bar{A}^1$            &          0.645272    &                   1                &  3.058        & 3223.1             &     9856.26 \\
		Brownian, $\bar{A}^2$                &        0.613716      &                 0.59717       &     3.31087    &    1274.01        &         4218.07  \\
		Brownian, $\bar{A}^3$              &      0.645702        &               0.959848       &    2.77357  &   3229.06        &        8956.03 \\
	\end{tabular}
	\caption{A summary of different statistics of the mM-MCMC method with $\varepsilon=10^{-6}$ for six combinations of the (approximate) macroscopic invariant distribution and macroscopic proposal moves. The different columns are the same as in Table~\ref{tab:threeatomgain1e-04}.}
	\label{tab:threeatomgain1e-06}
\end{table}

\paragraph{\textbf{Numerical results}}
The numerical results in Tables~\ref{tab:threeatomgain1e-04} and~\ref{tab:threeatomgain1e-06} indicate a large efficiency gain of mM-MCMC over the microscopic MALA algorithm. The efficiency gain is on the order of the time-scale separation $\sim \mathcal{O}(\varepsilon^{-1})$. For instance, a gain of a factor $8255$ in Table~\ref{tab:threeatomgain1e-06} indicates that the mM-MCMC method needs $8255$ times fewer sampling steps to obtain the same variance on the estimated mean of $\theta$ than the microscopic MALA method, for the same runtime. 

First, note that for both values of $\varepsilon$, the macroscopic acceptance rate is lower when using Brownian macroscopic proposals than when using Langevin dynamics proposals. This result is intuitive since the Brownian motion does not take into account the underlying macroscopic probability distribution, while the Langevin dynamics will automatically choose reaction coordinate values in regions of higher macroscopic probability. However, with a lower macroscopic acceptance rate comes a lower runtime as well since we need to reconstruct fewer microscopic samples and hence fewer evaluations of the microscopic potential energy. This effect is indeed visible in the fourth column where the runtime gain is higher for Brownian motion than that of the corresponding Langevin dynamics. Further, we note that there are large variance gains between the Langevin proposals and Brownian proposals for the second and third approximate free energy. We explains this as follows: the more exploration happens at the macroscopic levels, the more diverse the microscopic samples as well. As a result, expected values at the microscopic level are better represented and have a lower variance. Hence, the total efficiency gain of mM-MCMC is almost completely determined by the lower runtime due to the macroscopic Brownian proposals when the approximate distribution lies close to the true marginal distribution of the reaction coordinate. However, when there are significant differences between the approximate and exact macroscopic distributions, the choice of proposal moves does have a significant impact on the overall efficiency.

The choice of macroscopic invariant distribution $\bar{\mu}_0$, however, has a larger impact on the efficiency of mM-MCMC with direct reconstruction. First of all, one can see that the macroscopic acceptance rate is less affected by the choice of macroscopic invariant distribution than by the choice of macroscopic proposal move $q_0$. Indeed, the Langevin proposals are based on the free energy of their respective approximate macroscopic distribution. However, the microscopic acceptance rate is significantly affected by the approximate macroscopic distribution $\bar{\mu}_0$. Since the microscopic acceptance criterion is used to correct the microscopic samples from having the wrong macroscopic distribution, the more the approximate macroscopic distribution $\bar{\mu}_0$ deviates from the exact macroscopic distribution $\mu_0$, the lower the microscopic acceptance rate will be. Indeed, the microscopic acceptance rate for the approximate free energy $\bar{A}^1$ is much lower than that for $A$ since the local minima of $\bar{A}^1$ are located at different positions. On the other hand, the approximate free energy $\bar{A}^2$ lies closer to $A$ since only the height of both peaks different, resulting in a microscopic acceptance rate close to $1$.

Consequently, the closer $\bar{\mu}_0$ lies to the exact invariant distribution $\mu_0$ of the reaction coordinates, the higher the gain in variance is over the microscopic MALA method. Indeed, if the microscopic acceptance rate is low, we store the same microscopic sample many times, prohibiting a thorough exploration of the microscopic state space and thus keeping the variance obtained by mM-MCMC high. The choice of approximate macroscopic distribution $\bar{\mu}_0$ has a larger impact on the efficiency gain of mM-MCMC than the choice of macroscopic transition distribution $q_0$.

\subsubsection{Impact of $\bar{A}$ and $\bar{\nu}$ on the efficiency gain} \label{subsubsec:Anu}
\paragraph{\textbf{Experimental setup}}
For the third experiment on the three-atom molecule, we investigate the effect of the choice of reconstruction distribution $\bar{\nu}$ on the efficiency gain of mM-MCMC over MALA, in conjunction with the same three choices for the approximate free energy $\bar{A}$~\eqref{eq:freeenergies_2}. For consistency of the numerical results, we employ macroscopic proposal moves based on a gradient descent in $\theta$, with given approximate free energy $\bar{A}$. The two reconstruction distributions that we consider in this numerical experiment are
\begin{align} \label{eq:reconstructions}
\bar{\nu}^1(x|\theta) &= \nu(x|\theta) \nonumber \nonumber \\
\bar{\nu}^2(x|\theta) &\propto \exp\left(-\frac{(x_a-1)^2}{4\varepsilon}\right) \exp\left(-\frac{(r_c-1)^2}{4\varepsilon}\right).
\end{align}
The first reconstruction distribution is the exact time-invariant distribution~\eqref{eq:exactnu}, while the second distribution is obtained by increasing the variance on $x_a$ and $r_c$ by a factor of $2$, relative to $\bar{\nu}^1$.

\begin{table}[h]
	\centering
	\begin{tabular}{c|c|c|c|c|c}
		\centering
		Parameters & \pbox{15cm}{Macroscopic \\ acceptance rate} & \pbox{15cm}{Microscopic \\ acceptance rate} & \pbox{15cm}{Runtime \\ gain} & \pbox{15cm}{Variance \\ gain} & Total efficiency gain \\
		\hline
		$A, \bar{\nu}^1$ &                            0.749932        &               1        &           2.45692         &    85.3266         &        209.64 \\
		$\bar{A}^1, \bar{\nu}^1$                &            0.730384     &                  0.432508   &         2.62306           &  28.4122             &     74.527 \\
		$\bar{A}^2, \bar{\nu}^1$        &                    0.749653    &                   0.950238      &      1.95663            & 99.7186        &         195.112 \\
		$A, \bar{\nu}^2$ &                           0.749891       &                0.476963       &     2.45479            & 37.5736         &         92.2354 \\
		$\bar{A}^1, \bar{\nu}^2$ &                           0.730482                   &    0.266464        &    2.62922           &  17.2829              &    45.4405 \\
		$\bar{A}^2, \bar{\nu}^2$     &                       0.749654               &        0.474436       &     1.92062            & 50.257          &         96.5246 \\
	\end{tabular}
	\caption{A summary of different statistics of the mM-MCMC method with $\varepsilon=10^{-4}$ for six combinations of the (approximate) macroscopic invariant distribution and reconstruction distribution. For each of these six combinations, we record the average acceptance rate at the macroscopic level, the average acceptance rate on the microscopic level, conditioned on all accepted macroscopic samples and the gain in runtime and variance on the estimated mean of $\theta$ of mM-MCMC over the microscopic MALA method. The final column record the total efficiency gain of mM-MCMC, which is the product of the two former columns.}
	\label{tab:threeatomgainnu}
\end{table}

\begin{table}[h]
	\centering
	\begin{tabular}{c|c|c|c|c|c}
		\centering
		Parameters & \pbox{15cm}{Macroscopic \\ acceptance rate} & \pbox{15cm}{Microscopic \\ acceptance rate} & \pbox{15cm}{Runtime \\ gain} & \pbox{15cm}{Variance \\ gain} & Total efficiency gain \\
		\hline
		$A, \bar{\nu}^1$            &                0.749906          &             1             &      2.50343           & 3297.65          &         8255.44 \\
		$\bar{A}^1, \bar{\nu}^1$                 &           0.730538             &          0.43242    &         2.64621            & 933.64         &          2470.61\\
		$\bar{A}^2, \bar{\nu}^1$  &                          0.749596          &             0.950308        &    1.99301           & 2461.08             &      4904.96 \\
		$A, \bar{\nu}^2$              &              0.750029            &           0.47711    &         2.58486       &     1250.1  &                  3231.34 \\
		$\bar{A}^1, \bar{\nu}^2$        &                    0.730455         &              0.266443       &     2.65114         &    547.826           &       1452.37 \\
		$\bar{A}^2, \bar{\nu}^2$  &                          0.749564   &                    0.474496    &        1.93858     &       1396.05     &              2706.35 \\
	\end{tabular}
	\caption{A summary of different statistics of the mM-MCMC method with $\varepsilon=10^{-6}$ for six combinations of the (approximate) macroscopic invariant distribution and reconstruction distribution. The columns are the same as in Table~\ref{tab:threeatomgainnu}.}
	\label{tab:threeatomgainnu6}
\end{table}

In this experiment, we again compute the macroscopic acceptance rate, the microscopic acceptance rate, the gain in runtime for a fixed number of sampling steps, the gain in variance of the estimated mean of $\theta$ and the total efficiency gain of mM-MCMC over the microscopic MALA algorithm for each of the six combinations of the approximate macroscopic distribution $\bar{\mu}_0$ and reconstruction distribution $\bar{\nu}$. We perform each experiment with $N=10^6$ steps, the temperature parameter is $\beta=1$, the macroscopic time step is again $\Delta t=0.01$. We use a time step $\delta t = \varepsilon$ for the MALA algorithm. For a good comparison, the numerical results are averaged over $100$ independent runs. The numerical results are shown in Table~\ref{tab:threeatomgainnu} for $\varepsilon=10^{-4}$ and in Table~\ref{tab:threeatomgainnu6} for $\varepsilon=10^{-6}$.

\paragraph{\textbf{Numerical results}}
As intuitively expected, the choice of reconstruction distribution has a negligible impact on the macroscopic acceptance rate, but it does have a significant effect on the microscopic acceptance rate. For both values of $\varepsilon$ and for the three choices of $\bar{A}$, reconstruction distribution $\bar{\nu}^2$ has a lower microscopic acceptance rate than $\bar{\nu}^1$. Since the distribution $\bar{\nu}^2$ is not the time-invariant reconstruction distribution~\eqref{eq:exactnu}, the microscopic acceptance criterion needs to correct for this wrong reconstruction, lowering the average microscopic acceptance rate. As we also noted in the previous experiment, a lower microscopic acceptance rate keeps the variance on the estimated mean of $\theta$ high and hence the total efficiency gain of mM-MCMC over the microscopic MALA algorithm is small.

\subsubsection{Efficiency gain as a function of $\varepsilon$} \label{subsubsec:effgainepsdirect}
\begin{figure}
	\centering
	\includegraphics[width=0.65\linewidth]{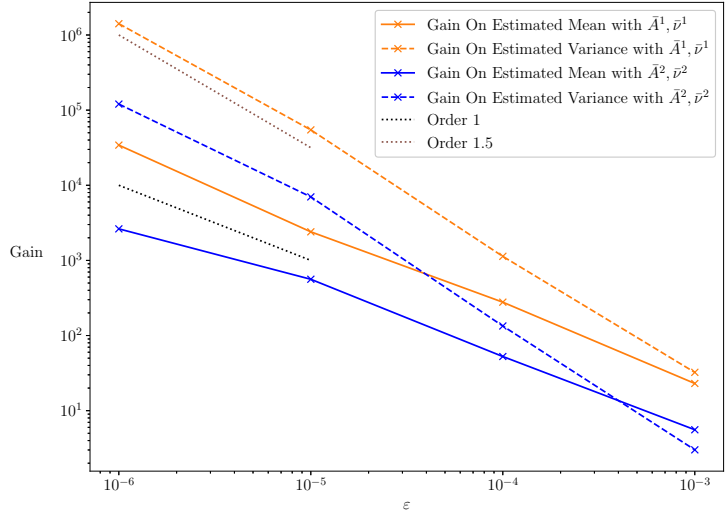}
	\caption{Efficiency gain of mM-MCMC over the standard MCMC method on the estimated mean (solid line) and variance (dashed line) of the angle $\theta$ for two different parameter choices: $(A, \bar{\nu}^1)$ (orange lines) and $(\bar{A}^1, \bar{\nu}^2)$ (blue lines). The gain is computed using the criterion[REF EFFICIENCY Crit], for $\varepsilon=10^{-i},  \ i=3,\dots,6$ and with $N=10^6$ microscopic samples. For both parameter choices, the efficiency gain on the estimated mean increases linearly with decreasing $\varepsilon$, and the gain on the estimated variance increases faster, approximately $\varepsilon^{-1.5}$. However, the efficiency gain on the estimated mean and variance of $\theta$ increases slower when the exact reconstruction property does not hold, i.e., $(\bar{A}^1, \bar{\nu}^2)$, than when it does, i.e., $(A, \bar{\nu}^1)$.}
	\label{fig:effgaindirect}
\end{figure}

\paragraph{\textbf{Experimental setup}}
In the fourth and final experiment on the three-atom molecule, we estimate the efficiency gain of mM-MCMC on the estimated mean and variance of the angle $\theta$, as a function of the time-scale separation $\varepsilon$. We consider two different choices for the approximate macroscopic distribution $\bar{\mu}_0$ and the reconstruction distribution $\bar{\nu}$. For the first choice, we take the exact free energy $A$[REF FREE energy 3-atom] and the exact time-invariant reconstruction distribution $\nu$ such that the exact reconstruction property holds. The other choice consists of approximate free energy $\bar{A}^1$ and reconstruction distribution $\bar{\nu}^2$. For both choices, the macroscopic proposal moves are based on the overdamped Langevin dynamics with the corresponding (approximate) free energy function and time step $\Delta t = 0.01$. We measure the efficiency gain for four values of the time-scale separation, $\varepsilon = 10^{-i}, \ i=3,\dots,6$ and with $N=10^6$ sampling steps. The microscopic time step for the microscopic MALA algorithm is $\delta t = \varepsilon$, and the inverse temperature is $\beta=1.$. The numerical results are shown in Figure~\ref{fig:effgaindirect}.

\paragraph{\textbf{Numerical results}}
For both choices of the parameters in the mM-MCMC method, the efficiency gain on the estimated mean and variance increases linearly or faster with decreasing $\varepsilon$, proving the mM-MCMC can accelerate the sampling of systems with a medium to large time-scale separation. Additionally, in case time-invariant $(A, \bar{\nu}^1)$, the efficiency gain is higher than when the exact reconstruction property does not hold, i.e., in case of $(\bar{A}^1, \bar{\nu}^2)$.

\subsection{Experiments with indirect reconstruction} \label{subsec:threeatomindirect}

\subsubsection{Impact of $\lambda$ on the efficiency of mM-MCMC} \label{subsubsec:threeatomlambda}
\begin{figure}
	\centering
	\includegraphics[width=0.85\linewidth]{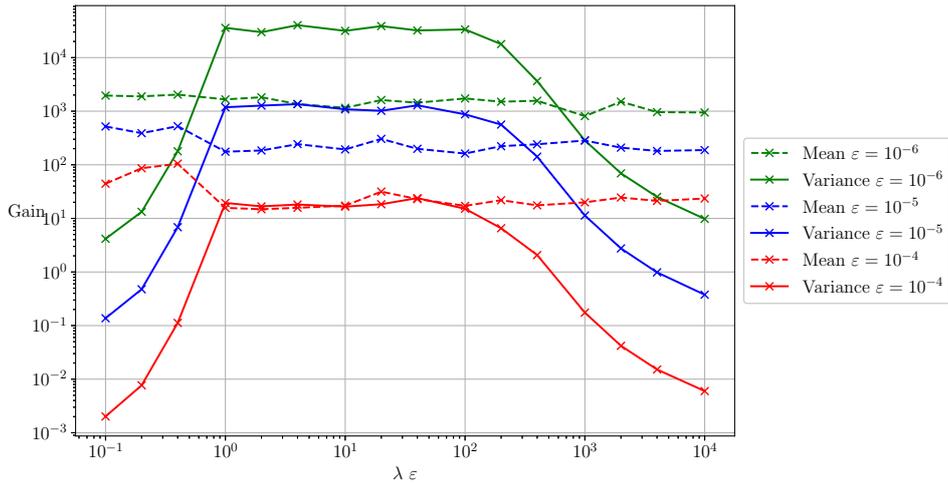}
	\caption{Efficiency gain of mM-MCMC on the estimate of the mean (full lines) and variance (dotted lines) of $\theta$, for different values of the time-scale separation. When $\lambda$ is larger $1/\varepsilon$ there is a clear efficiency gain. This gain is almost constant for a large range of values for $\lambda$. When $\lambda$ is too large, we need more biased steps to come near to the sampled value at the macroscopic level, and when $\lambda$ is smaller than $1/\varepsilon$, we take the sampled reaction value not enough into account and loose efficiency.}
	\label{fig:efflambda}
\end{figure}

\paragraph{Experimental setup}
It is paramount that the parameter $\lambda$ is chosen well, such that the mM-MCMC scheme is as efficiently as possible. He, we investigater the effect of the magnitude of $\lambda$ on the efficiency gain of the estimated mean and the estimated variance of $\theta$. Specifically, we compute the efficiency gain of mM-MCMC with indirect reconstruction over the microscopic MALA method for three values of the time-scale separation: $\varepsilon = 10^{-4}, \ 10^{-5}, \ \text{and} \  10^{-6}$. For a good comparison, we keep the number of biased time steps fixed at $K = 5$ and the step size of the biased simulation fixed at $\delta t=\min\{\varepsilon, \lambda^{-1}\}$ for stability. The other numerical parameters are $\beta=1$, the macroscopic time step is $\Delta t = 0.01$, the time step of the microscopic MALA method is $\varepsilon$ and we take $N=10^6$ sampling steps. Moreover, we compute the efficiency gain by averaging the estimated quantities over 100 independent runs. On Figure~\ref{fig:efflambda}, we show the efficiency gain of mM-MCMC as a function of $\lambda$, for multiple values of $\varepsilon$. Also, in Table~\ref{tab:threeatomgain_lambda}, we gather some statistics of the mM-MCMC method for the estimated variance of $\theta$ for several values of $\lambda$ and a fixed value of $\varepsilon=10^{-6}$.

\begin{table}[!h]
	\centering
	\begin{tabular}{c|c|c|c|c|c}
		\centering
		$\lambda  \cdot \varepsilon$ & \pbox{15cm}{Macroscopic \\ acceptance rate} & \pbox{15cm}{Microscopic \\ acceptance rate} & \pbox{15cm}{Runtime \\ gain} & \pbox{15cm}{Variance \\ gain} & \pbox{15cm}{Total \\ efficiency gain} \\
		\hline
		$0.1$ &                0.749749   &                    0.993677    &       0.218974         &   0.500816          &       0.109666\\
		$1$ & 0. 0.74966                     &   0.993468     &      0.233651  &     153735        &            35920.4 \\
		$10$ & 0.7496                         &0.993463      &     0.22258       & 141816            &        31565.4 \\
		$100$ & 0.749546                     &  0.993421     &      0.219963  &     152517           &         33548.1 \\
		$1000$ & 0.749553        &               0.993469    &       0.218501  &       1303.67       &            284.853 \\
	\end{tabular}
	\caption{Several statistics of the performance of mM-MCMC with indirect reconstruction over the microscopic MALA method when estimating the variance of $\theta$, for several values of $\lambda$. The time-scale separation parameter is fixed at $\varepsilon=10^{-6}$. The macroscopic and microscopic acceptance rates and the increase in execution time remain constant with varying $\lambda$. As is also visible from the greed solid line on Figure~\ref{fig:efflambda}, when $\lambda < \varepsilon^{-1}$, the efficiency gain is low because there is no gain on the variance at all. However, from the point $\lambda > \varepsilon^{-1}$, there is a significant variance reduction on the estimated variance of $\theta$, and therefore a significant increase in efficiency gain. Finally, when $\lambda$ is too large, relative to the fastest modes in the potential energy of the system, the total efficiency gain decreases again due to an increase of the variance on the estimated variance of $\theta$ by mM-MCMC.}
	\label{tab:threeatomgain_lambda}
\end{table}

\paragraph{Numerical results}
From previous sections, we know that efficiency gain increases when $\varepsilon$ decreases,. Second, when $\lambda < 1/\varepsilon$ there is no efficiency gain at all. Indeed, as we intuitively mentioned in Section~\ref{subsubsec:threeatomlambda}, when $\lambda$ is smaller than the stiffest mode in the system, the reaction coordinate in the biased dynamics will be less driven towards the value sampled at the macroscopic level. We are then effectively ignoring the macroscopic MCMC step so that there is no variance reduction on expectations of $\theta$. This effect is also visible in the first row of Table~\ref{tab:threeatomgain_lambda} where there is indeed no gain on the variance on the estimated variance of $\theta$. On the other hand, when $\lambda$ is very large (larger than $10^2 / \varepsilon$ in this case), the efficiency also starts to decrease as we are not simulating enough biased steps for the reaction coordinate value to approximate the sampled value at the macroscopic level well. This effect is especially visible on the efficiency gain of the estimate variance of $\theta$ in the bottom row of Table~\ref{tab:threeatomgain_lambda}. We would therefore need more than $5$ biased steps to equilibrate around each sampled reaction coordinate value at the macroscopic level of the mM-MCMC algorithm, also reducing the efficiency. To conclude, there is a large range of values for $\lambda$ that give a large and almost identical efficiency gain (middle rows of Table~\ref{tab:threeatomgain_lambda}). This range, between $1/\varepsilon$ and $100/\varepsilon$, is the same for a large range of $\varepsilon$ values. In practice, it is therefore a good idea to choose $\lambda$ approximately on the order of the stiffest mode of the molecular system.

\subsubsection{Efficiency gain as a function of time-scale separation} \label{subsubsec:triatomindriecteps}

\begin{table}
	\centering
	\begin{tabular}{c|c|c|c|c|c}
		\centering
		$\varepsilon$ & \pbox{15cm}{Macroscopic \\ acceptance rate} & \pbox{15cm}{Microscopic \\ acceptance rate} & \pbox{15cm}{Runtime \\ gain} & \pbox{15cm}{Variance \\ gain} & \pbox{15cm}{Total \\ efficiency gain} \\
		\hline
		$10^{-3}$ & 0.749975           &            0.993528      &     0.212558       &      10.3927         &         2.20905\\
		$10^{-4}$ & 0.74936          &              0.993588     &      0.211415      &       68.6401          &       14.5115 \\
		$10^{-5}$ & 0.750197          &             0.993299    &       0.212561     &       920.651        &         195.695 \\
		$10^{-6}$ & 0.749498              &         0.993405  &         0.237227   &        7041.69          &       1670.48 \\
	\end{tabular}
	\caption{A summary of different statistics of the mM-MCMC method when applied to estimating the mean of $\theta$, for multiple values of $\varepsilon$. We record the macroscopic and microscopic acceptance rates, the runtime and variance gains of the mM-MCMC method with indirect reconstruction and the total efficiency gain. Note that the cost of pre-computing the free energy is not included in the runtime.}
	\label{tab:threeatomgain_mean}
\end{table}

\begin{figure}
	\centering
	\includegraphics[width=0.55\linewidth]{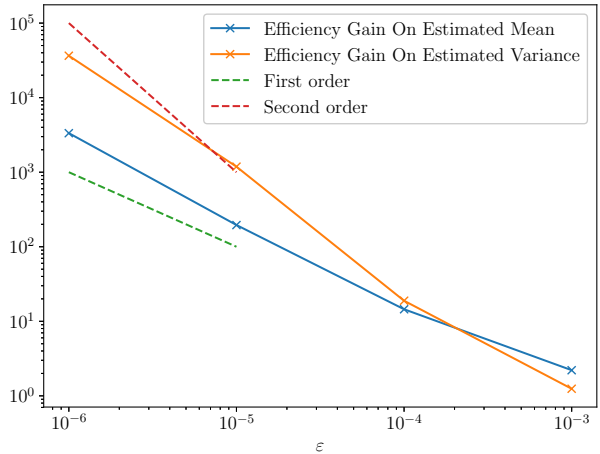}
	\caption{Efficiency gain of mM-MCMC over standard MCMC as a function of the time-scale separation. The efficiency gain for the mean of $\theta$ (blue) increases linearly with decreasing $\varepsilon$, while the efficiency gain for the variance of $\theta$ increases faster than linearly, although slower than quadratically.}
	\label{fig:gainepstriatom}
\end{figure}

\paragraph{Experimental setup}
Following the previous experiment, we numerically compare the efficiency gain of the mM-MCMC method over the microscopic MALA method for different values of $\varepsilon$. We let the small-scale parameter $\varepsilon$ vary between $10^{-6}$ and $10^{-3}$ and we consider the efficiency gain on the estimated mean and the estimated variance of $\theta$. The efficiency gains on the estimated mean and variance of $\theta$ are depicted on Figure~\ref{fig:gainepstriatom} as a function of $\varepsilon$. We also display the averaged macroscopic and microscopic acceptance rates, the reduction in runtime, as well as the reduction of the variance on the estimated quantities of interest obtained by mM-MCMC with indirect reconstruction in Table~\ref{tab:threeatomgain_mean} where we consider the estimated mean of $\theta$, and in Table~\ref{tab:threeatomgain_variance} for the estimated variance of $\theta$.

\paragraph{Numerical results}
First, on Figure~\ref{fig:gainepstriatom}, we see that the efficiency gain of mM-MCMC increases at least linearly with decreasing $\varepsilon$. For large values of $\varepsilon$, one can see that there is almost no efficiency gain at all. The reason for this behaviour is clear from Tables~\ref{tab:threeatomgain_mean} and~\ref{tab:threeatomgain_variance}. From the third column of these tables, we conclude that the runtime of mM-MCMC with indirect reconstruction is larger than that of the microscopic MALA method due to the computational overhead of the biased simulation. Hence, for large values of $\varepsilon$, this computational overhead (third column) is dominant over the reduction in variance (fourth column). However, when $\varepsilon$ is small, the reduction in variance by mM-MCMC is dominant over the increase in computational time. Currently, we do not have analytic expressions for the gain as a function of the time-scale separation, but Figure~\ref{fig:gainepstriatom} clearly shows the merit of mM-MCMC with indirect reconstruction for medium and large time-scale separations.

\begin{table}
	\centering
	\begin{tabular}{c|c|c|c|c|c}
		\centering
		$\varepsilon$ & \pbox{15cm}{Macroscopic \\ acceptance rate} & \pbox{15cm}{Microscopic \\ acceptance rate} & \pbox{15cm}{Runtime \\ gain} & \pbox{15cm}{Variance \\ gain} & \pbox{15cm}{Total \\ efficiency gain} \\
		\hline
		$10^{-3}$ & 0.749975              &         0.993528  &         0.212558     &        5.85671         &         1.24489\\
		$10^{-4}$ & 0.74936     &                  0.993588      &     0.211415     &       89.2989                &  18.8791 \\
		$10^{-5}$ & 0.750197         &              0.993299        &   0.212561    &      5580.75      &            1186.25 \\
		$10^{-6}$ & 0.749498          &             0.993405      &     0.237227  &      153706           &         36463.2 \\
	\end{tabular}
	\caption{A summary of different statistics of the mM-MCMC method when applied to estimating the variance of $\theta$, for multiple values of $\varepsilon$. The conclusions on the macroscopic and microscopic acceptance rates and the increase in execution time are the same as in Table~\ref{tab:threeatomgain_mean}.}
	\label{tab:threeatomgain_variance}
\end{table}

As a second observation, note that the macroscopic acceptance rates (first column of both tables) is independent of $\varepsilon$, as we intuitively may expect. Also, we see that the microscopic acceptance rate after indirect reconstruction (second column) is close to $1$ so that only little redundant computational work is performed during indirect reconstruction. This result shows that the indirect reconstruction is an efficient technique to reconstruct a microscopic sample close to a sub-manifold of constant reaction coordinate value.

\subsubsection{Comparison of mM-MCMC with direct and indirect reconstruction}\label{subsubsec:comparisondirectindirect}

\begin{figure}
	\centering
	\includegraphics[width=0.55\linewidth]{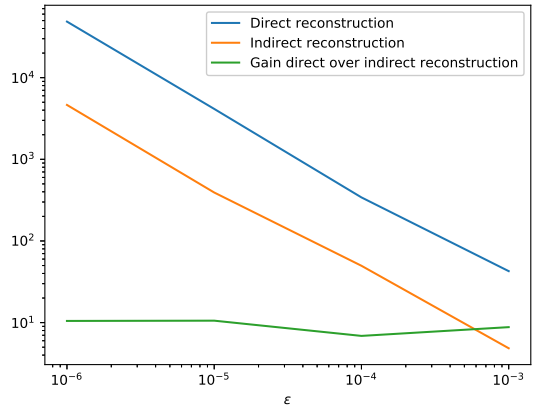}
	\caption{Efficiency gain of mM-MCMC with direct reconstruction (blue) and indirect reconstruction (orange) over the microscopic MALA method as a function of $\varepsilon$. The green curve measures the efficiency gain of the direct reconstruction algorithm over the indirect reconstruction method. One can see that the efficiency gain is a constant factor of less than $10$ lower then the gain made by the direct reconstruction variant.}
	\label{fig:directindirect}
\end{figure}

Having studied the performance of mM-MCMC with indirect reconstruction for multiple values of the time-scale separation, we now compare the performance of this method to its direct reconstruction variant. We expect the direct reconstruction algorithm to be faster for a given number of sampling steps due to the computational overhead of the biased dynamics, while the reduction in variance on estimated quantities should be almost identical. 

\paragraph{Experimental setup}
In Figure~\ref{fig:directindirect}, we depict the efficiency gain of both mM-MCMC variants over the microscopic MALA method on the estimated mean of $\theta$ for a large range of values of $\varepsilon$. In addition, we also plot the efficiency gain of the direct reconstruction algorithm over \emph{the indirect reconstruction algorithm} for $N=10^6$ sampling steps. We fix the macroscopic time step at $\Delta t = 0.02$ for both mM-MCMC variants for a good comparison. The numerical parameters for the indirect reconstruction method are $K = 5$, $\lambda = \varepsilon^{-1}$ and $\delta t = \varepsilon$, and the time step of the microscopic MALA method is also $\varepsilon$. For both variants, we use the exact marginal distribution $\mu_0$ at the macroscopic level and we use Brownian proposal moves. We also use reconstruction distribution $\nu$ for the direct reconstruction method. In Table~\ref{tab:directindirect}, we also show the gain in runtime, the gain in variance on the estimated mean of $\theta$ and the total efficiency gain of mM-MCMC with direct reconstruction over mM-MCMC with indirect reconstruction.

\paragraph{Numerical results}
The efficiency gain of mM-MCMC with indirect reconstruction is a constant factor lower than the efficiency gain of its direct reconstruction variant, independent of the time-scale separation. If we diagnose this effect more carefully in Table~\ref{tab:directindirect}, one can see that the lower efficiency gain is purely due to the larger runtime of mM-MCMC with indirect reconstruction. The decrease in variance on the estimated mean of $\theta$ is almost the same.

\begin{table}
	\centering
	\begin{tabular}{c|c|c|c}
		\centering
		$\varepsilon$  & \pbox{15cm}{Runtime \\ gain} & \pbox{15cm}{Variance \\ gain} & \pbox{15cm}{Total \\ efficiency gain} \\
		\hline
		$10^{-3}$  &         10.0329     &     0.877206        &    8.80089\\
		$10^{-4}$  &   10.1001     &     0.680705     &       6.87522 \\
		$10^{-5}$  &   10.1374     &     1.04204   &         10.5635 \\
		$10^{-6}$  &    9.99052 &        1.57429     &       10.4853 \\
	\end{tabular}
	\caption{A summary of different statistics that summarize the efficiency gain of mM-MCMC with direct reconstruction over the indirect reconstruction algorithm. The efficiency gain of the direct reconstruction algorithm is almost the same for a large range of values of $\varepsilon$ and this efficiency gain is completely due to the lower runtime of the direct reconstruction. Both variants obtain the exact same variance reduction over the microscopic MALA method.}
	\label{tab:directindirect}
\end{table}

\section{Butane} \label{subsec:butane}
\begin{figure}
	\centering
	\includegraphics[width=0.45\linewidth]{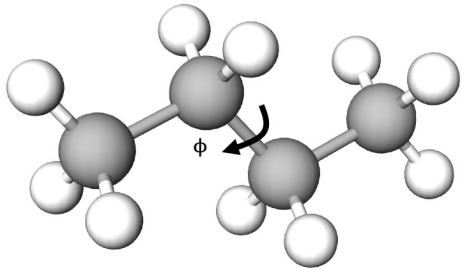}
	\caption{The butane molecule. The carbon atoms are grey and hydrogen is white.}
	\label{fig:butane}
\end{figure}

We now test the mM-MCMC method with direct reconstruction on a somewhat larger method to show that its overall efficiency gain extends to other molecules. We specifically focus on the behaviour of the torsion angle $\phi$ (see Figure~\ref{fig:butane}), as this angle determines the global conformation of the whole molecule.

The full potential energy of butane reads
\begin{equation} \label{eq:Vbutane}
V(x) =  \frac{1}{2} \sum_{i\in bonds} k_b \ (r_i(x) - r_0)^2 + \frac{1}{2} \sum_{i \in angles} k_a \ (\theta_i(x) - \theta_0)^2 + V_{\phi}(\phi(x)),
\end{equation}
with parameteres $k_ b =  1.17 \cdot 10^{6}$, $k_a = 62500$,  $r_0 = 1.53 $, $\theta_0 = 112^{\circ}$ and potential energy of the torsion angle beign
\[
V_{\phi}(\phi) = c_0 + c_1\cos(\phi) + c_2 \cos(\phi)^2 + c_3 \cos(\phi)^3.
\]
The remaining paramters have values $c_0=1031.36 ,  c_1 = 2037.82, c_2 = 158.52, c_3 = -3227.7 $. The time-invariant distribution of butane is then a result of~\eqref{eq:gibbs}.

Given the slow nature of the torsion angle and the faster vibrations of the bonds and regular angles, it makes sense to choose our reaction coordinate to be the torsion angle, i.e., 
\[
\xi(x) = \phi.
\]
The free energy of this reaction coordinate is readily visible from the potential energy of the full system~\eqref{eq:Vbutane}
\begin{equation} \label{eq:freeenergy_butane}
A(\phi) =V_{\phi}(\phi),
\end{equation}
and the time-invariant reconstruction distribution is then given by~\eqref{eq:exactnu}.

\paragraph{Outline of this section}
We adapt a couple of experiments that we performed on the three-atom molecule to the problem of butane. We already performed a visual comparison between the microscopic MALA and the mM-MCMC methods in our previous work~\cite{vandecasteele2022}. There, we demonstrated that the mM-MCMC method with direct reconstructoin can achieve significant efficiency gains on butane. We take this analysis a step further, and compare the macroscopic reaction coordinate sampler (MALA) and the mM-MCMC method in case there is a significant bias on the approximate macroscopic distribution in Section~\ref{subsubsec:macro_butane}. Afterwards, we compare the efficiency gains of the mM-MCMC method for different choices of $\bar{\mu}_0$ and $q_0$, Section\ref{subsubsec:impact_butane_mu_q}, and different choices of $\bar{\mu}_0$ and $\bar{\nu}$, Section~\ref{subsubsec:impact_butane_mu_nu}.

\subsection{Comparison with the macroscopic reaction coordinate sampler} \label{subsubsec:macro_butane}
As a second numerical illustration of the merit of mM-MCMC, we compare the sampling results between a macroscopic sampler and the mM-MCMC method. We have already discovered that the mM-MCMC method can correct for large biases in the approximate macroscopic distribution in case of the three-atom molecule. The butane test case is, however, more challenging because the local minima of $A(z)$ are farther apart and two of the local minima have a lower weight than the third one.

\paragraph{Experimental setup}
We define two approximate free energy functions that induce approximate macroscopic distributions
\begin{equation} \label{eq:approx_free_butane}
\begin{aligned}
\bar{A}^1(\phi) &= c_0 + c_1\cos(1.2\phi) + c_2 \cos(1.2\phi)^2 + c_3 \cos(1.2\phi)^3, \\
\bar{A}^2(\phi) &= c_0 + c_1\cos(\phi) + c_2 \cos(\phi)^2 + c_3 \cos(\phi)^3  + 500\cos(\phi-1).
\end{aligned}
\end{equation}
The first approximate free energy contracts the local minima of the exact free energy~\eqref{eq:freeenergy_butane} by a factor $1.2^{-1}$, while the third approximate free energy gives a larger relative weight to most left local minimum, and decreases the weight of the other peaks.

\begin{figure}
	\centering
	\begin{subfigure}[b]{0.5\textwidth}
		\centering
		\includegraphics[width=1.\linewidth]{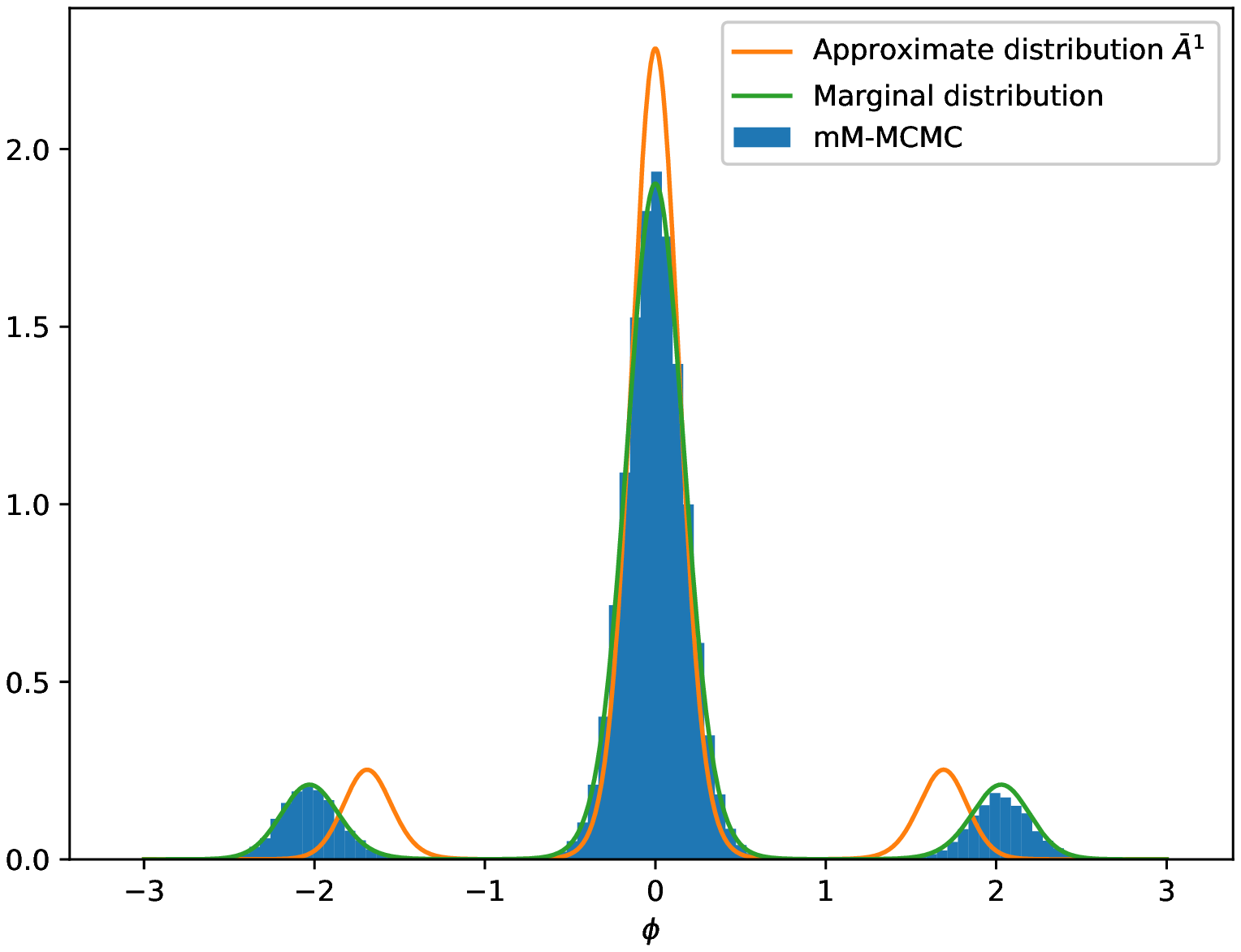}
	\end{subfigure}%
	\begin{subfigure}[b]{0.5\textwidth}
		\centering
		\includegraphics[width=1.\linewidth]{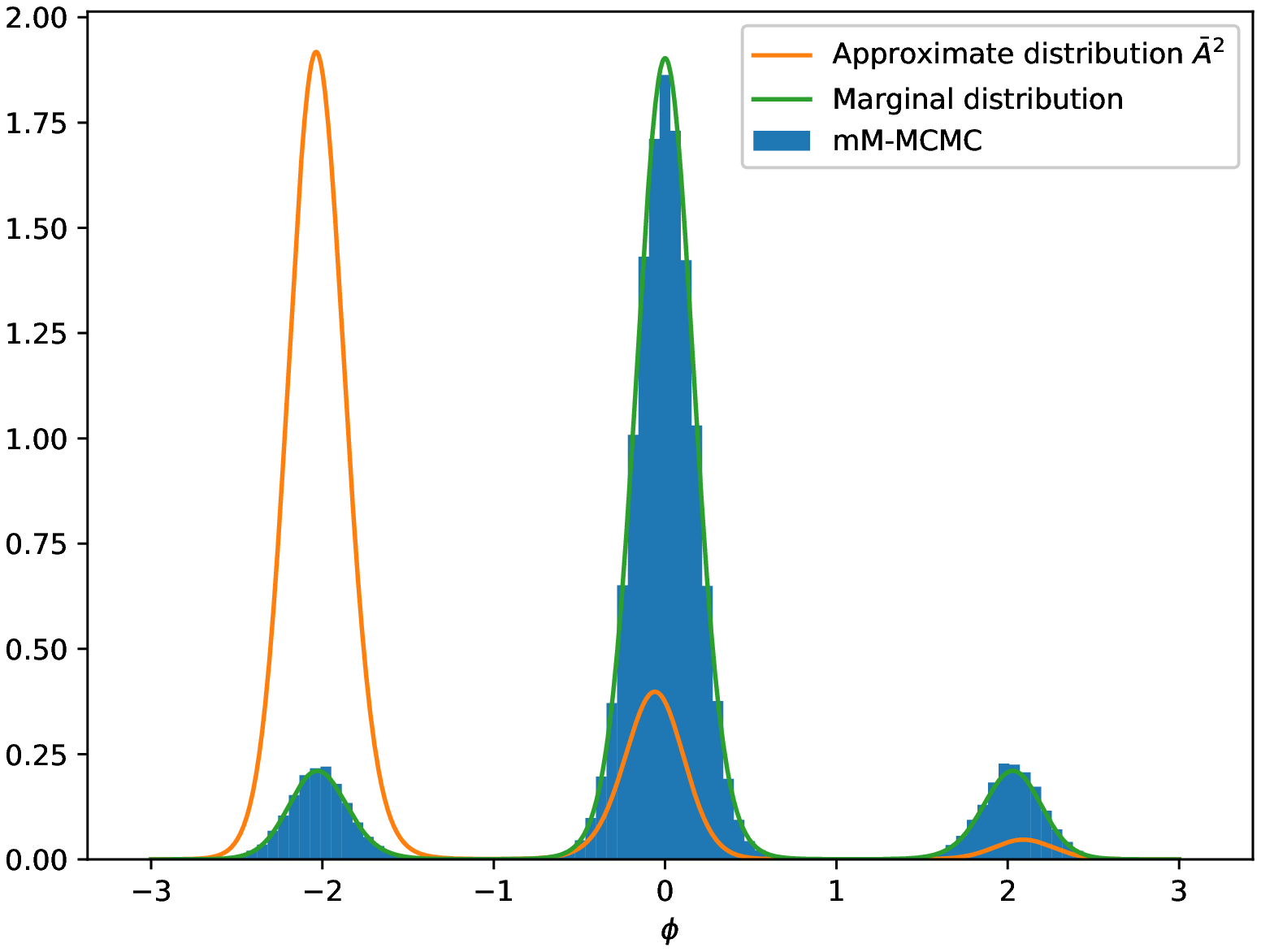}
	\end{subfigure}
	\caption{Histogram of the mM-MCMC method (blue) for two different approximate macroscopic distributions: $\bar{A}^1$ (left) and $\bar{A}^2$ (right). The orange line depicts the approximate macroscopic distribution used in each experiment, and the green line represents the marginal distribution of the reaction coordinate values. We can see that the mM-MCMC is able to overcome a bad approximate macroscopic distribution and samples the marginal distribution correctly.}
	\label{fig:macro_approx_butane}
\end{figure}

We now investigate the numerical performance of mM-MCMC over a macroscopic sampler for each of the above defined approximate distribution. On the macroscopic level, we employ a MALA method with time step $\Delta t=5\cdot 10^{-4}$ to generate proposals for both the macroscopic and mM-MCMC samplers. We also use the time-invariant reconstruction distribution $\nu$ to generate microscopic samples for mM-MCMC. The numerical results for approximate free energy $\bar{A}^1$ and $\bar{A}^2$ are shown in Figure~\ref{fig:macro_approx_butane}.

\paragraph{Numerical results}
Clearly, the micro-macro MCMC method is also able to overcome bad approximate macroscopic distributions for butane. This effect is especially visible on the right figure. Indeed, the approximate macroscopic distribution gives a completely wrong idea about the marginal reaction coordinate distribution by assigning a too large weight to the left-most peak, and a too small weight to the middle peak. However, the micro-macro MCMC corrects for this misrepresentation in the microscopic accept/reject step.

\subsection{Impact of $\bar{A}$ and $q_0$ on the efficiency gain} \label{subsubsec:impact_butane_mu_q}
We have shown in the previous section that the mM-MCMC can correct a bad approximate macroscopic distribution and still sample the marginal distribution of the reaction coordinate $\phi$. In this section, we take this experiment one step further by comparing the actual efficiency gains of six combinations of the approximate macroscopic distribution $\bar{\mu}_0$ and the macroscopic proposals moves $q_0$. To measure the efficiency of mM-MCMC over the microscopic MALA scheme, we numerically estimate the variance of $\phi$ with the obtained microscopic samples. We fix the reconstruction distribution to the exact marginal $\nu$ in this experiment. We consider three approximate macroscopic distributions in this experiment: the exact free energy~\eqref{eq:freeenergy_butane} and the two approximate free energies defined in the previous section~\eqref{eq:approx_free_butane}. The two macroscopic proposal moves are defined by a discretization of the dynamics
\begin{equation*}
\begin{aligned}
q_0^1 &: d\phi = -\nabla \bar{A}(\phi) dt + \sqrt{2 \beta^{-1}} dW \\
q_0^2 &: d\phi = \sqrt{2 \beta^{-1}} dW,
\end{aligned}
\end{equation*}
where $\bar{A}$ is a shorthand for each of the three (approximate) macroscopic free energies used in this experiment. It can easily be seen that the first proposal move represents the MALA method~\eqref{eq:mala}, while the second proposal move is a simple scaled Brownian motion.

\paragraph{Experimental setup}
For a good comparison, we average the estimated means of $\phi$ over $100$ independent runs In each independent mM-MCMC and MALA run we use $N=10^6$ sampling steps, the macroscopic time step for mM-MCMC is $\Delta t=5 \ 10^{-4}$ and the microscopic time step for the MALA method is $\delta t=10^{-6}$. We also keep track of the macroscopic and microscopic acceptance rates of the mM-MCMC method, the runtime of both methods for a given number of sampling steps, the variance on the estimated variance of $\phi$ for both method and the total efficiency gain of mM-MCMC over MALA. These numerical results are displayed in Table~\ref{tab:butane_Aq}. 

\begin{table}[h]
	\centering
	\begin{tabular}{c|c|c|c|c|c}
		\centering
		Parameters & \pbox{15cm}{Macroscopic \\ acceptance rate} & \pbox{15cm}{Microscopic \\ acceptance rate} & \pbox{15cm}{Runtime \\ gain} & \pbox{15cm}{Variance \\ gain} & \pbox{15cm}{Total \\ efficiency gain} \\
		\hline
		$q_0^1$, $A$               & 0.289293  &  1             &  23.8795  & 1050.6     &  25087.9\\
		$q_0^2$, $A$    & 0.423083  &  1             &  19.4699  &  5156.4    &  100396 \\
		$q_0^1$, $\bar{A}^1$    & 0.313243  &  0.297524  &  22.8156 &  0.085769 &  1.95688 \\
		$q_0^2$, $\bar{A}^1$               & 0.407327  &  0.754024  & 19.6211  &  0.364752 &   7.15682 \\
		$q_0^1$ $\bar{A}^2$     & 0.30793    &  0.762832  &  21.7298 & 1341.02    &  29140\\
		$q_0^2$ $\bar{A}^2$     & 0.432885  &  0.822621  &  18.3769 &  3279.37   &   60264.5 \\
	\end{tabular}
	\caption{The macroscopic and microscopic acceptance rates, the runtime and variance gain and the total efficiency gain of the mM-MCMC method for six combinations of the approximate macroscopic distribution and the macroscopic proposal moves. The closer the approximate macroscopic distribution to the true marginal of the reaction coordinate, the higher the efficiency gain. Macroscopic proposal moves based on Brownian increments also result in a larger efficiency gain than proposal moves based on the MALA method.}
	\label{tab:butane_Aq}
\end{table}

\paragraph{Numerical results}
It is clear from the numerical results that approximate macroscopic distribution with $\bar{A}^1$ results in a poorer efficiency gain than $\bar{A}^2$. The reason for this behaviour is that the outer peaks of the distribution associated to $\bar{A}^1$  are closer to $0$ than the actual local minima of the exact free energy $A$. At the macroscopic level, reaction coordinate samples are thus 'punished' for lying further from the outer most peaks, result in a low microscopic acceptance rate and hence efficiency gain. This problem does not happen with approximate free energy $\bar{A}^2$ as its outer peaks have the same location as $A$.

Second, similarly to the three-atom molecule, macroscopic proposal moves based on Brownian increments result in a larger efficiency gain than those based on the macroscopic MALA method. Indeed, although the runtime gain is lower for Brownian proposals due to a larger macroscopic acceptance rate, the variance gain is factor $3$ to $5$ higher for $q_0^2$ than for $q_0^1$.

\subsection{Impact of $\bar{A}$, and $\bar{\nu}$ on the efficiency gain} \label{subsubsec:impact_butane_mu_nu}
Finally, we investigate the role of the reconstruction distribution $\bar{\nu}$ on the efficiency gain of the mM-MCMC method. To this end, we investigate the overall performance of mM-MCMC for all six combinations consisting of three choices for the (approximate) free energy as in the previous section, and two choices for the reconstruction distribution:
\begin{equation*}
\begin{aligned}
\bar{\nu}^1(x|z) &= \nu(x|z) = \mu(x) \mu_0(z)^{-1} \\
\bar{\nu}^2(x|z) &= \mu(x)^{\frac{1}{2}}  Z(z)^{-1}.
\end{aligned}
\end{equation*}
Here, the distribution $\bar{\nu}^2$ is obtained by doubling the variances of each of the microscopic degrees of freedom of butane, and $Z(z)$ is the normalization constant.

\paragraph{Experimental setup}
For each of the above combinations, we compute the efficiency gain on the variance of $\phi$. The macroscopic time step for mM-MCMC is $\Delta t = 5 \cdot 10^{-4}$ and the microscopic time step for the MALA method is $\delta t=10^{-6}$. We use $N=10^6$ sampling steps for both method and we average each combination over $100$ independent runs. The macroscopic and microscopic acceptance rates, the runtime and variance gains and the total efficiency gain of mM-MCMC over the microscopic MALA  method are shown in Table~\ref{tab:butane_Anu}.

\begin{table}[h]
	\centering
	\begin{tabular}{c|c|c|c|c|c}
		\centering
		Parameters & \pbox{15cm}{Macroscopic \\ acceptance rate} & \pbox{15cm}{Microscopic \\ acceptance rate} & \pbox{15cm}{Runtime \\ gain} & \pbox{15cm}{Variance \\ gain} & \pbox{15cm}{Total \\ efficiency gain} \\
		\hline
		$A$, $\bar{\nu}^1$    & 0.423083  &  1             &  19.4699  &  5156.4    &  100396 \\
		$A$, $\bar{\nu}^2$    &0.422962 & 0.154385 &  9.88807    & 0.90294    &    8.92833\\
		$\bar{A}^1$, $\bar{\nu}^1$              & 0.407327  &  0.754024  & 19.6211  &  0.364752 &   7.15682 \\
		$\bar{A}^1$, $\bar{\nu}^2$              &0.408502  & 0.130944 &  10.1383   &   1.09534    &   11.1049 \\
		$\bar{A}^2$, $\bar{\nu}^1$    & 0.432885  &  0.822621  &  18.3769 &  3279.37   &   60264.5 \\
		$\bar{A}^2$, $\bar{\nu}^2$    &0.432576  & 0.14719  &  9.74243   &  2.5751   &     25.0878  \\
	\end{tabular}
	\caption{The macroscopic and microscopic acceptance rates, the runtime and variance gain and the total efficiency gain of the mM-MCMC method for six combinations of the approximate macroscopic distribution and the reconstruction distribution. Reconstruction distribution $\bar{\nu}^2$, obtained by enlarging the variance of each of the microscopic components of butane by a factor $2$, significantly affects the total efficiency gain.}
	\label{tab:butane_Anu}
\end{table}

\paragraph{Numerical results}
The choice of reconstruction distribution has a significant effect on the total efficiency gain of the mM-MCMC method. Indeed, for both macroscopic distributions given by the (approximate) free energy functions $A$ and $\bar{A}^2$, the efficiency gain is much larger in case of $\bar{\nu}^1$. This effect is not visible for approximate free energy $\bar{A}^1$, but the efficiency gains are too small to be interesting. The numerical results indicate the need for an efficient and general way of reconstruction microscopic samples on (or close to) the sub-manifold of a given reaction coordinate value. We explore such a method in the next section.

\subsection{Indirect reconstruction: optimal $\lambda$} \label{subsec:optimal_lambda}
Before testing the mM-MCMC scheme with indirect reconstruction, we need to determine a good value of $\lambda = 1/\varepsilon$ first. From~\eqref{eq:Vbutane}, $k_b$ is  the largest force constant,  so we will base $\lambda$  on that.

\paragraph{Experimental setup} We expect $\lambda$ to be of  the same order of  $k_b$. We  therefore test  a  range  of  values  between  $0.1 kb$  to $100 k_b$. For  each  such  value  of  $\lambda$, we do one hundred independend  runs of the mM-MCMC algorithm, and compute the average efficiency gain~\ref{eq:effgain}. We measure both the efficiency  gain on the mean torsion angle and variance of the torsion angle. We also run the microscopic sampler a hundred times. The parameters of the  mM-MCMC method are $N=10^5$ and $K=10$.  The macroscopic sampler uses Brownian increments with variance $\Delta t=0.4$, and the reconstruction time step is $\delta t = 0.4/\lambda$. The microscopic sampler uses  $10^5$ steps with a  time step of $1.0/k_b$ for stability.

\begin{figure}
	\centering
	\includegraphics[width=0.7\linewidth]{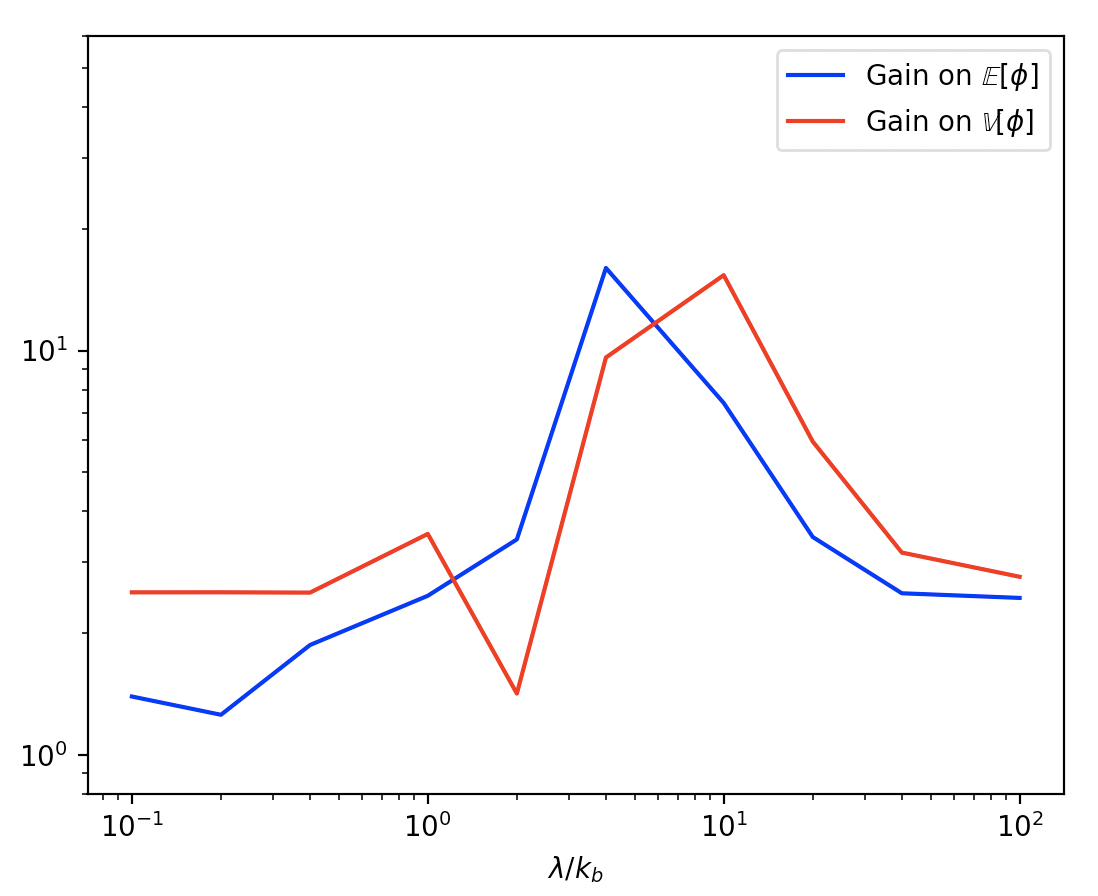}
	\caption{Efficiency gain of mM-MCMC over the microscopic MALA method as a function of $\lambda$.  Blue: efficiency gain on the mean torsion angle. Red: efficiency gain on the variance of the torsion angle.}
	\label{fig:lambdaselectionbutane}
\end{figure}

\paragraph{Numerical results}
The results are shown in Figure~\ref{fig:lambdaselectionbutane}. For all values we tried here, the mM-MCMC method with indirect reconstruction always obtains an efficiency gain. Note that this gain is not necesarily maximal because we did not tune the reconstruction time step $\delta t$ optimally. For the interpretation of the results, there is a clear area between $4k_b$ and $10k_b$ where both the efficiency gains on the mean and variance of $\phi$ are maximised. Any of these $\lambda$-values can be used during samping. Note that this result is similar to that in section~\ref{subsubsec:threeatomlambda} in that the best values of $\lambda$ are indeed on the order of  $k_b$, but larger values are also possible.

\subsection{Comparison between macroscopic and microscopic samplers}
In this final experiment, we compare the mM-MCMC method to both the microscopic sampler and the macroscopic sampler. We already showed that mM-MCMC is more efficient than the microscopic sampler, and the goal of this experiment is to prove that mM-MCMC is as accurate as a macroscopic sampler. Note that it is impossible to use a macroscopic sampler in most applications because the free energy is not directly available.

\paragraph{Experimental setup}
We run all three algorithms for  $N=10^5$  sampling steps. The macroscopic sampler uses Brownian increments with variance $\Delta t = 0.4$. As in previous experiments, the microscopic sampler uses MALA steps with variance $\delta t = 2/k_b$. Next, the mM-MCMC method also uses macroscopic  Brownian increments with $\Delta  t=0.4$ and MAlA reconstruction steps with $\delta t = 0.4/\lambda$. From previous figure, we choose $\lambda = 4 k_b$ and $K = 10$.  For each of these algorithms, we plot the Gaussian kernel density estimation (KDE) with bandwidth $h=0.03$. The bandwidth was chosen so that it captures the right peaks of the distribution, while not adding superfluous detail due to the finite nature of the samplers.

\begin{figure}
	\centering
	\includegraphics[width=0.7\linewidth]{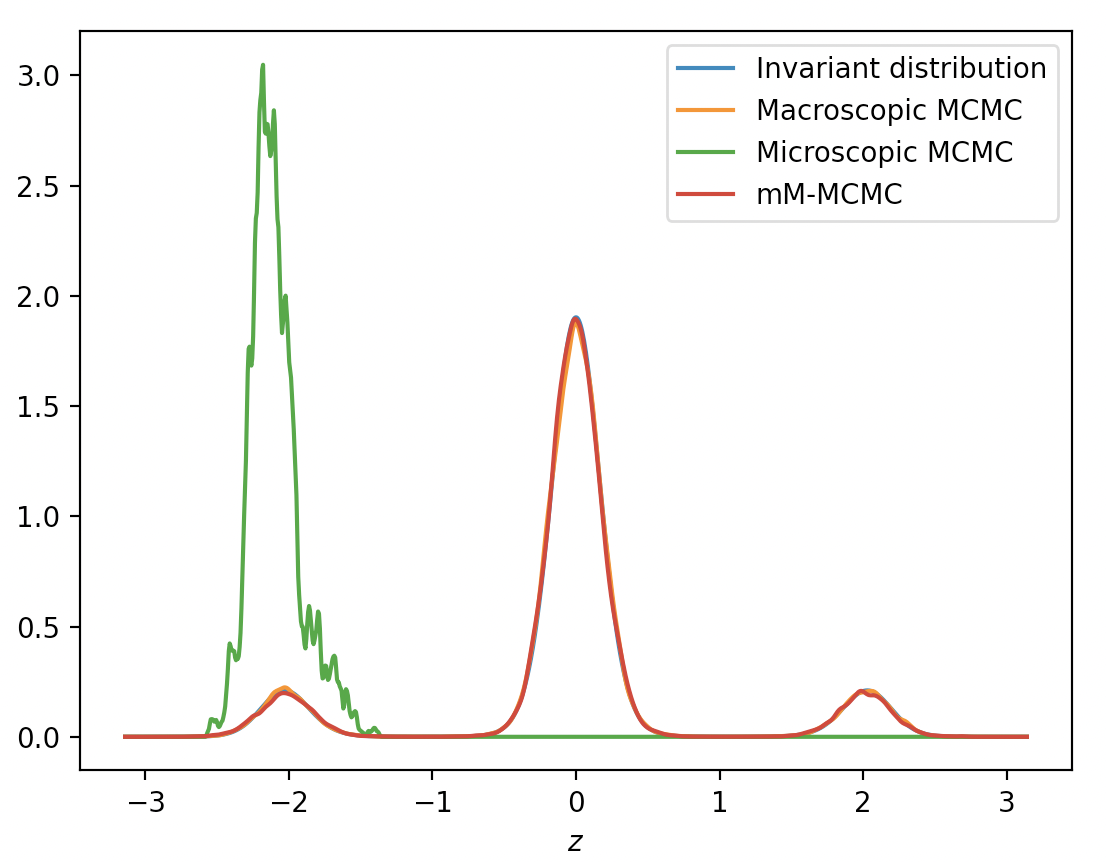}
	\caption{Gaussian kernel density estimation with bandwidth $h=0.03$ for: the macroscopic MCMC method (orange), the microscopic MCMC method (green), and the mM-MCMC method with indirect reconstruction (red). The true free energy is drawn in blue.}
	\label{fig:kdem}
\end{figure}

\paragraph{Numerical results}
Figure~\ref{fig:kdem} shows the Gaussian KDE's for all three algorithms: Macroscopic (orange), microscopic(green) and mM-MCMC(red). The true invariant distribution of the torsion angle is shown in blue. On the figure, we see that both distributions of the macroscopic sampler and mM-MCMC method are very close to invariant distribution, while the microscopic sampler is stuck in the left lobe. We conclude that the mM-MCMC method indeed samples the correct invariant distribution, and obtains a large efficiency gain over the incorrect microscopic sampler.

\section{Conclusion} \label{sec:conclusion}
We examined the effect of various paramters of the mM-MCMC scheme on the overall efficiency gain of said method. Particularly, in case of direct reconstruction, we investigated how $\bar{A}$, $q_0$ and $\bar{\nu}$ changed the efficiency on two molecules: the three-atom molecule and butane. The conclusion is identical in both cases. The total efficiency gain is largest when we choose $\bar{A}$ close to the exact free energy, $\bar{\nu}$ close to $\nu$ and when the macroscopic proposals are based on Brownian motion (no gradient information). The first choice may be difficult to accomplish in practice because of the high cost associated to computing the free energy. However, an approximate distribution that has a similar structure works well too. We also analyzed the mM-MCMC with indirect reconstruction. First, we determined the optimal value of the biasing potential, $\lambda$ and found that the optimal value is related to the fastest mode in the system. Second, we also saw that the mM-MCMC scheme with indirect reconstruction is as good as the macroscopic sampler, but with the added benefit of having microscopic samples available. Finally, a comparision between both mM-MCMC variants showed that the direct variant is faster, per sample, but both have the same variance gain, though indirect reconstruction is applicable in every application.

\paragraph{Acknowledgment}
This work was supported by the Flemish Fund for Scientific Research (Grant 1179820N).

\bibliographystyle{plain}
\bibliography{bib}

\begin{thebibliography}{10}

\bibitem{berg1991multicanonical}
B.~A. Berg and T.~Neuhaus.
\newblock Multicanonical algorithms for first order phase transitions.
\newblock {\em Physics Letters B}, 267(2):249--253, 1991.

\bibitem{christen2005markov}
A.~J. Christen and C;. Fox.
\newblock Markov chain monte carlo using an approximation.
\newblock {\em Journal of Computational and Graphical statistics},
  14(4):795--810, 2005.

\bibitem{comer2015adaptive}
J.~Comer, J.~C. Gumbart, J.~H{\'e}nin, T.~Leli{\`e}vre, A.~Pohorille, and
  C.~Chipot.
\newblock The adaptive biasing force method: Everything you always wanted to
  know but were afraid to ask.
\newblock {\em The Journal of Physical Chemistry B}, 119(3):1129--1151, 2015.

\bibitem{darve2001calculating}
E.~Darve and A.~Pohorille.
\newblock Calculating free energies using average force.
\newblock {\em The Journal of chemical physics}, 115(20):9169--9183, 2001.

\bibitem{darve2008adaptive}
E.~Darve, D.~Rodr{\'\i}guez-G{\'o}mez, and A.~Pohorille.
\newblock Adaptive biasing force method for scalar and vector free energy
  calculations.
\newblock {\em The Journal of chemical physics}, 128(14):144120, 2008.

\bibitem{dickson2010free}
B.~M. Dickson, F.~Legoll, T.~Leli{\`e}vre, G.~Stoltz, and P.~Fleurat-Lessard.
\newblock Free energy calculations: An efficient adaptive biasing potential
  method.
\newblock {\em The Journal of Physical Chemistry B}, 114(17):5823--5830, 2010.

\bibitem{dodwell2019multilevel}
Tim~J Dodwell, Christian Ketelsen, Robert Scheichl, and Aretha~L Teckentrup.
\newblock Multilevel markov chain monte carlo.
\newblock {\em Siam Review}, 61(3):509--545, 2019.

\bibitem{hastings1970monte}
Keith~W. Hastings.
\newblock Monte carlo sampling methods using markov chains and their
  applications.
\newblock 1970.

\bibitem{jasra2018markov}
A.~Jasra, K.~Law, and Y.~Xu.
\newblock Markov chain simulation for multilevel monte carlo.
\newblock {\em arXiv preprint arXiv:1806.09754}, 2018.

\bibitem{kalligiannaki2012coupled}
E.~Kalligiannaki, M.~A. Katsoulakis, and P.~Plech{\'a}{\v{c}}.
\newblock Coupled coarse graining and markov chain monte carlo for lattice
  systems.
\newblock In {\em Numerical Analysis of Multiscale Computations}, pages
  235--257. Springer, 2012.

\bibitem{laio2002escaping}
A.~Laio and M.~Parrinello.
\newblock Escaping free-energy minima.
\newblock {\em Proceedings of the National Academy of Sciences},
  99(20):12562--12566, 2002.

\bibitem{le2012mathematical}
Claude Le~Bris, Tony Leli{\`e}vre, Mitchell Luskin, and Danny Perez.
\newblock A mathematical formalization of the parallel replica dynamics.
\newblock 2012.

\bibitem{legoll2012some}
Fr{\'e}d{\'e}ric Legoll and Tony Leli{\`e}vre.
\newblock Some remarks on free energy and coarse-graining.
\newblock In {\em Numerical Analysis of Multiscale Computations}, pages
  279--329. Springer, 2012.

\bibitem{marsili2006self}
S.~Marsili, A.~Barducci, R.~Chelli, P.~Procacci, and V.~Schettino.
\newblock Self-healing umbrella sampling: a non-equilibrium approach for
  quantitative free energy calculations.
\newblock {\em The Journal of Physical Chemistry B}, 110(29):14011--14013,
  2006.

\bibitem{metropolis1953equation}
Nicholas Metropolis, Arianna~W. Rosenbluth, Marshall~N. Rosenbluth, Augusta~H.
  Teller, and Edward Teller.
\newblock Equation of state calculations by fast computing machines.
\newblock {\em The journal of chemical physics}, 21(6):1087--1092, 1953.

\bibitem{vandecasteele2022}
Hannes Vandecasteele and Giovanni Samaey.
\newblock A micro-macro markov chain monte carlo method for molecular dynamics
  using reaction coordinate proposals.
\newblock 2022.

\bibitem{wang2001efficient}
F.~Wang and D.~P. Landau.
\newblock Efficient, multiple-range random walk algorithm to calculate the
  density of states.
\newblock {\em Physical review letters}, 86(10):2050, 2001.

\bibitem{xifara2014langevin}
Tatiana Xifara, Chris Sherlock, Samuel Livingstone, Simon Byrne, and Mark
  Girolami.
\newblock Langevin diffusions and the metropolis-adjusted langevin algorithm.
\newblock {\em Statistics \& Probability Letters}, 91:14--19, 2014.

\bibitem{zheng2008random}
L.~Zheng, M.~Chen, and W.~Yang.
\newblock Random walk in orthogonal space to achieve efficient free-energy
  simulation of complex systems.
\newblock {\em Proceedings of the National Academy of Sciences},
  105(51):20227--20232, 2008.

\end{thebibliography}

\end{document}